\documentclass{article}    
\usepackage{latexsym}
\usepackage{amsbsy}
\usepackage{amssymb}
\usepackage[latin1]{inputenc} 
\usepackage{amscd}

\begin{document}                                                                       

\newcommand{\ep}{\hspace*{\fill}$\Box$}
\newcommand{\eps}{\varepsilon}
\newcommand{\pr}{{\it Proof. }}
\newcommand{\ms}{\medskip\\}
\newcommand{\cl}{\mbox{\rm cl}}
\newcommand{\R}{\mathbb R}
\newcommand{\N}{\mathbb N}
\newcommand{\C}{\mathbb C}
\newcommand{\Z}{\mathbb Z}
\newcommand{\K}{\mathbb K}
\newcommand{\sR}{\mathbb R}
\newcommand{\sN}{\mathbb N}
\newcommand{\gs}{\ensuremath{{\mathcal G}} }
\newcommand{\gso}{\ensuremath{{\mathcal G}(\Omega)} }
\newcommand{\gsrn}{\ensuremath{{\mathcal G}(\R^n)} }
\newcommand{\gsrp}{\ensuremath{{\mathcal G}(\R^p)} }
\newcommand{\gsrq}{\ensuremath{{\mathcal G}(\R^q)} }
\newcommand{\gst}{\ensuremath{{\mathcal G}_\tau} }
\newcommand{\gsto}{\ensuremath{{\mathcal G}_\tau(\Omega)} }
\newcommand{\gstrn}{\ensuremath{{\mathcal G}_\tau(\R^n)} }
\newcommand{\gstrp}{\ensuremath{{\mathcal G}_\tau(\R^p)} }
\newcommand{\gstrq}{\ensuremath{{\mathcal G}_\tau(\R^q)} }
\newcommand{\es}{\ensuremath{{\mathcal E}} }
\newcommand{\esm}{\ensuremath{{\mathcal E}_M} }
\newcommand{\ns}{\ensuremath{{\mathcal N}} }
\newcommand{\nso}{\ensuremath{{\mathcal N}(\Omega)} }
\newcommand{\est}{\ensuremath{{\mathcal E}_\tau} }
\newcommand{\nst}{\ensuremath{{\mathcal N}_\tau} }
\newcommand{\ks}{\ensuremath{{\mathcal K}} }
\newcommand{\rs}{\ensuremath{{\mathcal R}} }
\newcommand{\cs}{\ensuremath{{\mathcal C}} }
\newcommand{\comp}{\subset\subset}
\newcommand{\cinfty}{{\cal C}^\infty}
\newcommand{\catensor}{\ensuremath{{\mathcal A}^\otimes} }
\newtheorem{lemma}{Lemma}
\newcommand{\blem}{\begin{lemma}}
\newcommand{\elem}{\end{lemma}}
\newtheorem{theorem}{Theorem}
\newcommand{\btheorem}{\begin{theorem}}
\newcommand{\etheorem}{\end{theorem}}
\newtheorem{prop}{Proposition}
\newcommand{\bp}{\begin{prop}}
\newcommand{\eprop}{\end{prop}}
\newtheorem{cor}{Corollary}
\newcommand{\bc}{\begin{cor}}
\newcommand{\ecor}{\end{cor}}
\newtheorem{example}{Example}
\newcommand{\bex}{\begin{example} \rm }
\newcommand{\eex}{\end{example}}
\newtheorem{defi}{Definition}
\newcommand{\bd}{\begin{defi}}
\newcommand{\edefi}{\end{defi}}
\newtheorem{remark}{Remark}
\newcommand{\brem}{\begin{remark} \rm}
\newcommand{\erem}{\end{remark}}

\newcommand{\beast}{\begin{eqnarray*}}
\newcommand{\eeast}{\end{eqnarray*}}
\newcommand{\wsc}[1]{\overline{#1}^{wsc}}
\newcommand{\todo}[1]{$\clubsuit$\ {\tt #1}\ $\clubsuit$}
\newcommand{\rem}[1]{\vadjust{\rlap{\kern\hsize\thinspace\vbox%
                       to0pt{\hbox{${}_\clubsuit${\small\tt #1}}\vss}}}}
\newcommand{\ahat}{\ensuremath{\hat{\mathcal{A}}_0(M)} }
\newcommand{\atil}{\ensuremath{\tilde{\mathcal{A}}_0(M)} }
\newcommand{\aqtil}{\ensuremath{\tilde{\mathcal{A}}_q(M)} } 
\newcommand{\ehat}{\ensuremath{\hat{\mathcal{E}}(M)} } 
\newcommand{\emhat}{\ensuremath{\hat{\mathcal{E}}_m(M)} }
\newcommand{\nhat}{\ensuremath{\hat{\mathcal{N}}(M)} }
\newcommand{\ghat}{\ensuremath{\hat{\mathcal{G}}(M)} } 
\newcommand{\lhat}{\ensuremath{\hat{L}_X} }                    
\newcommand{\SSS}{{\cal S}}
\newcommand{\al}{\alpha}
\newcommand{\bet}{\beta} 
\newcommand{\ga}{\gamma}
\newcommand{\Om}{\Omega}\newcommand{\Ga}{\Gamma}\newcommand{\om}{\omega}
\newcommand{\si}{\sigma}\newcommand{\la}{\lambda}
\newcommand{\de}{\delta}
\newcommand{\vphi}{\varphi}\newcommand{\dl}{{\displaystyle \lim_{\eta>0}}\,}
\newcommand{\intl}{\int\limits}\newcommand{\su}{\sum\limits_{i=1}^2}
\newcommand{\D}{{\cal D}}\newcommand{\Vol}{\mbox{Vol\,}}
\newcommand{\Or}{\mbox{Or}}\newcommand{\sign}{\mbox{sign}}
\newcommand{\na}{\nabla}\newcommand{\pa}{\partial}
\newcommand{\ti}{\tilde}\newcommand{\T}{{\cal T}} \newcommand{\G}{{\cal G}}
\newcommand{\DD}{{\cal D}}\newcommand{\X}{{\cal X}}\newcommand{\E}{{\cal E}} 
\newcommand{\CC}{{\cal C}}\newcommand{\vo}{\Vol}
\newcommand{\bat}{\bar t}
\newcommand{\bx}{\bar x}
\newcommand{\by}{\bar y} \newcommand{\bz}{\bar z}\newcommand{\br}{\bar r}
\newcommand{\fr}{\frac{1}}\newcommand{\il}{\int\limits}
\newcommand{\nn}{\nonumber}
\newcommand{\supp}{\mathop{\mathrm{supp}}}
\newcommand{\vp}{\mbox{vp}\frac{1}{x}}\newcommand{\A}{{\cal A}}
\newcommand{\Ll}{L_{\mbox{\small loc}}}\newcommand{\Hl}{H_{\mbox{\small loc}}}
\newcommand{\Lll}{L_{\mbox{\scriptsize loc}}}
\newcommand{\be}{ \begin{equation} }\newcommand{\ee}{\end{equation} }
\newcommand{\beq}{ \begin{equation} }\newcommand{\eeq}{\end{equation} }
\newcommand{\bea}{\begin{eqnarray}}\newcommand{\eea}{\end{eqnarray}}
\newcommand{\beas}{\begin{eqnarray*}}\newcommand{\eeas}{\end{eqnarray*}}
\newcommand{\beqs}{\begin{equation*}}\newcommand{\eeqs}{\end{equation*}}
\newcommand{\lb}{\label}\newcommand{\rf}{\ref}
\newcommand{\GL}{\mbox{GL}}\newcommand{\bfs}{\boldsymbol}
\newcommand{\ben}{\begin{enumerate}}\newcommand{\een}{\end{enumerate}}
\newcommand{\ba}{\begin{array}}\newcommand{\ea}{\end{array}}
\newcommand{\id}{\mathop{\mathrm{id}}}
\newcommand{\lgl}{\langle}
\newcommand{\rgl}{\rangle}
\newcommand{\ca}{{\cal A}}
\newcommand{\cc}{{\cal C}}
\newcommand{\cd}{{\cal D}}
\newcommand{\cg}{{\cal G}}
\title{Foundations of a Nonlinear Distributional Geometry} 
\author{Michael Kunzinger\footnote{ e-mail: michael.kunzinger@univie.ac.at}, 
Roland Steinbauer\footnote{ e-mail: roland.steinbauer@univie.ac.at}\\
Department of Mathematics, University of Vienna\\
Strudlhofg. 4, A-1090 Wien, Austria.}
\date{September, 2001}

\maketitle

\begin{abstract}
Colombeau's construction of generalized functions (in its special variant) is
extended to a theory of generalized sections of vector bundles. As particular 
cases, generalized tensor analysis and exterior algebra are studied. A point value
characterization for generalized functions on manifolds is derived, several
algebraic characterizations of spaces of  generalized sections are established and consistency
properties with respect to linear distributional geometry are derived. An application
to nonsmooth mechanics indicates the additional flexibility offered by this approach
compared to the purely distributional picture.
\vskip1em
\noindent{\footnotesize {\bf Mathematics Subject Classification (2000):} 
Primary: 46F30; secondary: 46T30, 46F10, 37J05}

\noindent{\footnotesize{\bf Keywords:} Algebras of generalized functions, Colombeau algebras, 
generalized sections of vector bundles, distributional geometry }
\end{abstract}


\section{Introduction}\label{intro}
After their introduction in \cite{c1}, \cite{c2} the main applications of Colombeau's new generalized
functions lay in the field of linear and nonlinear partial differential equations involving
singular coefficients or data (cf.\ \cite{MObook}, \cite{esiproc} and the literature
cited therein for a survey). Over the past few years, however, the theory has found 
a growing number of applications in a more geometric context, most notably in general
relativity (cf.\ e.g., \cite{clarke}, \cite{ultra}, \cite{herbertgeo},\cite{penrose},  as well as
\cite{vickersESI} for a survey). This shift of focus has necessitated a certain restructuring
of the fundamental building blocks of the theory in order to adapt to the additional 
requirement of diffeomorphism invariance. Only recently (\cite{found}, \cite{vim}) this
task has been completed for the scalar case. To be precise, this restructuring took
place in the framework of the so-called {\em full} Colombeau algebra, 
distinguished by the existence of a canonical embedding of the space of Schwartz 
distributions into the algebra.

Already at a very early stage of development the
so-called {\em special}  (or simplified) variant of Colombeau's algebras was introduced 
(cf.\ e.g., \cite{AB}). This variant of the construction does not allow for a canonical 
embedding of the space of distributions. If 
such an embedding is needed, then the full version of the theory (\cite{vim}) 
should be employed. On the other hand, in the special version due to its simpler basic structure 
(elements are basically equivalence classes of nets of smooth functions) an adaptation of
geometric constructions from the smooth setting to the generalized functions framework 
can be carried out more directly than in the full variant. In particular,
diffeomorphism invariance of the basic building blocks of the construction is 
automatically satisfied. Moreover, in applications where a distinguished regularization process is
available (e.g., due to certain symmetries of the problem under consideration or due
to a smoothing procedure suggested by physics) it is
often preferable to work in the special setting. 
Consequently there has been an increasing number of applications of the special algebra
to geometric problems (cf.\ e.g., \cite{RD}, \cite{DPS}, \cite{symm}, \cite{penrose}).
The aim of the present paper is to initiate a systematic development of global analysis
in this setting. 

For an alternative approach to algebras of generalized functions on manifolds
based on the abstract differential geometry developed by A.\ Mallios (\cite{mallios}) 
we refer to \cite{mallios+r}. 

The plan of the present work is as follows: In the remainder of the current section we 
fix some notation concerning differential geometry, while in section~\ref{specol} we 
recall the basic facts on special Colombeau algebras.
In section \ref{distgeom} we give a quick overview of distributional
geometry, introducing those constructions that later on will furnish our main objects 
of reference for the limiting behavior of the corresponding Colombeau objects. 
In section \ref{specalgmf} we introduce several equivalent definitions of as well 
as some basic operations on the special algebra of generalized functions $\G(X)$ 
on a manifold $X$. We then derive a point value characterization of elements of $\G(X)$,
a feature which distinguishes the present framework from the purely distributional one 
and serves as an important tool for generalizing notions from classical geometry.
Section \ref{assoz1} is devoted to a study of the compatibility of the current
approach with respect to distributional and smooth geometry. We discuss in detail 
the question of embedding $\D'(X)$ into $\G(X)$ reaching
the conclusion that a canonical and geometric embedding (in a sense to be
made precise there) indeed is not feasible. On the other hand we give a
simple construction of a (non-canonical) embedding that
extends to an injective sheaf morphism $\D'(\_)\hookrightarrow\G(\_)$ which
coincides with the natural (``constant'') embedding on $\CC^\infty(\_)$.
Furthermore we set up coupled calculus, in particular
the notion of $k$-association which is stronger than the notion of association used
in the local theory and---in the absence of a geometric embedding of 
$\D'$---serves to make precise statements on the compatibility with respect
to the distributional and $\CC^k$-setting. 
In section \ref{gsvb} we introduce generalized sections of vector bundles.
We prove some algebraic characterizations of these sheaves of $\G(\_)$-modules and
again establish consistency results with respect to the classical setting. Important special
cases of these general constructions are worked out in sections \ref{gta} (generalized
tensor analysis) and \ref{extalg} (exterior algebra). In particular, section \ref{extalg} 
provides an application to nonsmooth mechanics.
 
Notations from differential geometry will basically be chosen in accordance 
with \cite{AM}, \cite{michorbook}. Throughout this paper, $X$ will denote a
paracompact, smooth Hausdorff manifold of dimension $n$. For any vector bundle 
$E\to X$, by $\Gamma^k(X,E)$ (resp.\ $\Gamma^k_c(X,E)$) $(0\leq k\leq\infty)$
we denote the $\CC^k(X)$-module of (compactly supported) $\CC^k$-sections in $E$ and 
frequently drop the superscript if $k=\infty$. 
In particular, by $\mathfrak{X}(X)$ resp.\ $\Om^k(X)$ we denote the space
of smooth vector fields resp.\ $k$-forms on $X$.  
Generally, for $M_1,\dots,M_k$, $M_0$ modules over a commutative ring $R$, 
$L_R(M_1,\dots,M_k;M_0)$ denotes the $R$-module of $R$-$k$-linear maps from
$M_1\times\dots\times M_k$ into $M_0$. Since we will be considering tensor
products with respect to different rings $R$, the notation $M_1\otimes_R M_2$ 
will be used. By ${\cal P}(X,E)$ we denote the space of linear differential 
operators $\Gamma(X,E) \to \Gamma(X,E)$. For $E=X\times \R$ we write ${\cal P}(X)$
for ${\cal P}(X,E)$. 

\section{Special Colombeau algebras}\label{specol}

In this section we shortly recall some basic facts on algebras of generalized functions
and, in particular, Colombeau's so-called special construction on open sets of
Euclidean space. The key idea in constructing these algebras (which contain the space of 
Schwartz distributions and provide maximal consistency with respect to classical analysis) 
is regularization by nets of smooth functions and the use of asymptotic estimates with 
respect to the regularization parameter $\eps$. More precisely we employ a quotient 
construction as follows (for details we refer to \cite{AB}, \cite{c2}): 
denoting by $\Om$ an open subset of $\R^n$ we set (with $I=(0,1]$)
\begin{eqnarray*}
        {\cal E}(\Omega)&:=&(C^\infty(\Omega))^I\nn\\
        {\cal E}_M(\Omega)&:=&\{(u_\varepsilon)_{\varepsilon\in I}\in {\cal E}(\Omega)
                : \forall K\subset\subset\Omega, \forall\alpha\in\N_0^n
                \ \exists N\in \N \\
                &&\hspace{3cm}
                \sup_{x\in K}|\partial^\alpha u_\varepsilon
                (x)|=O(\varepsilon^{-N})\mbox{ as }\varepsilon\rightarrow 0\}\\
        {\cal N}(\Omega)&:=&\{(u_\varepsilon)_{\varepsilon\in I}\in {\cal E}(\Omega)
                : \forall K\subset\subset\Omega,\forall\al\in \N_0^n, \forall m\in\N:\\
                &&\hspace{3.3cm}
                \sup_{x\in K}|\partial^\alpha u_\varepsilon
                (x)|=O(\varepsilon^{m})\mbox{ as }\varepsilon\rightarrow 0\}.
\end{eqnarray*}
The {\em special Colombeau algebra on $\Om$} is defined as the quotient space
\[
        {\cal G}(\Omega):={\cal E}_M(\Omega)\,/\,{\cal N}(\Omega)\,.
\]
Since we will only be considering this type of algebras we will omit
the term ``special'' henceforth.
Elements of $\gso$ will be denoted by capital letters, representatives by small 
letters, i.e., $\gso\ni U=\cl[(u_\eps)_\eps]=(u_\eps)_\eps+\nso$. $\gs(\_)$ is a {\em fine
sheaf of differential algebras} containing the smooth functions on $\Omega$ as a
subalgebra embedded simply by $\sigma(f)=\cl[(f)_\eps]$. 

To embed non-smooth
distributions we first have to fix a mollifier $\rho\in\SSS(\R^n)$ with unit integral
satisfying the moment conditions $\int\rho(x)\,x^\al\,dx\,=\,0\ \forall 
|\al|\geq 1$. Setting $\rho_\eps(x)=(1/\eps^n)\rho(x/\eps)$, compactly supported distributions 
are embedded by $\iota_0(w)=((w*\rho_\eps|_\Om)_\eps+\nso$. Using partitions of unity and
suitable cut-off functions one may explicitly construct an embedding $\iota:\ \D'(\Om\hookrightarrow
\gso$ which naturally induces a unique sheaf morphism (of complex vector spaces) $\hat\iota:\D'(\_)
\hookrightarrow\gs(\_)$ extending $\iota_0$, commuting with partial derivatives and
its restriction to $\cinfty(\_)$ being a sheaf morphism of algebras.  
Note that $\hat\iota$ depends on the choice of the mollifier $\rho$, hence is 
non-canonical. This fact reflects a fundamental property of nonlinear modeling: 
In general, nonlinear properties of a singular object depend on the regularization. 
Additional input on the regularization from, say, a
physical model may enter the mathematical theory via this interface, leading to a 
sensible description of the problem at hand; in many cases this will be more natural
than the use of a ``canonical'' embedding of $\D'$ into $\G$.
 
A ``macroscopic'' description of calculations in $\G$ can often be effected
through the concept of {\em association:} $U$, $V\in \G(\Om)$ are called associated
if $u_\eps - v_\eps \to 0$ in $\D'(\Om)$, $U$ is called associated to $w\in \D'(\Om)$
if $u_\eps \to w$ in $\D'(\Om)$. Clearly these notions do not depend on the particular
representatives  and the first one gives rise to a linear quotient space of $\gso$,
which extends the notion of distributional equality to the level of the algebra. 

Finally, we note that inserting $x\in\Om$
into $U\in \G(\Om)$ componentwise yields a well-defined element of the ring of 
generalized numbers ${\cal K}$ (corresponding to $\K=\R$ resp.\ $\C$), defined
as the set of moderate nets of numbers ($(r_\eps)_\eps \in \K^I$ with 
$|r_\eps| = O(\eps^{-N})$ for some $N$) modulo negligible nets 
($|r_\eps| = O(\eps^{m})$ for each $m$).

\section{Distributional geometry} \label{distgeom}

We shortly recall the basic facts of distributional geometry, i.e., of the theory of
distribution valued sections of vector bundles.

On open sets of $\R^n$ a distribution is defined to be a continuous
linear functional on the (LF)-space of smooth, compactly supported
test {\em functions } $\vphi$. Any smooth (even any locally integrable)
{\em function} $f$ gives rise to a regular distribution via the
(natural) assignment $\vphi\mapsto \int f(x)\cdot\vphi(x)\,dx.$
On a general manifold $X$, these two statements cannot hold
simultaneously in a meaningful way (with emphasis on ``functions'').
In the absence of a preferred  measure the objects to be
integrated  are (one-){\em densities} which are sections of the
volume bundle $\Vol(X)$ (cf.\ e.g., \cite{simanca}). 
Thus, either the nature of test ``functions'' $\vphi$ or
of regular distributions $f$ or of both has to be changed in such a way
that their product $f\cdot\vphi$ becomes a density. 
Since the product of a density with a smooth function is again a density
there immediately arise two (in a sense, complementary) ways of proceeding.
On one hand, we can replace test {\it functions} by
test {\it densities} and define a distribution to be a continuous linear
functional on the space of these densities. Then again each (say, smooth)
{\it function} can be considered as a distribution.
This is in accordance with e.g., \cite{hoer1}, Sec.\ 6.3.
On the other hand, we could keep the function character of the test objects; 
then the regular objects in the dual space of the space of test {\em functions} 
have to be taken as (smooth) {\em densities }on $X$. This is the definition adopted e.g.,
in \cite{dieudonne3}, Ch.\ XVII. 
 
More generally, the burden of rendering $f\cdot\vphi$ a
density can be split up in one part contributed by $f$ and in
a (complementary) part contributed by $\vphi$. This is done by defining
for each real $q$ the notion of a $q$-density as a section of
the $q$-volume bundle $\Vol^q(X)$ of $X$ (see e.g., \cite{simanca,parker}). 
Moreover for arbitrary real $q,q'$, the product of a $q$-density with a $q'$-density
is a $(q+q')$-density; one-densities are just densities in the
above sense and zero-densities correspond to functions.
If we now define the test objects to be (compactly supported, smooth) 
$q$-densities, the appropriate $(1-q)$-densities can be embedded in their dual space
as regular objects. Note that the case $q=1/2$ is of particular interest due to the fact
that the product together with the integral induces a natural Hilbert space structure. 

The goal of defining vector valued distributions of a certain density character
finally is achieved by considering  $q$-densities with values in some vector bundle $E$ over $X$ 
as test objects, that is sections of the bundle $E\otimes\Vol^q(X)$.
An appropriate regular dual object for such (compactly supported, smooth) sections $u$
obviously would be a smooth section $f$ of the bundle $E^*\otimes\vo^{1-q}(X)$
where $E^*$ denotes the dual bundle of $E$; the canonical bilinear form $(\,.\,|\,.\,)$ 
on $E^*\times E$ and the product of densities make $(f|u )$ a one-density.
Interchanging $E$ and $E^*$ as well as $q$ and $1-q$, we finally arrive at the definition
of $E$-valued distributions of density character $q$ and order $k$ (to be formally given below)
as the dual of the space of compactly supported $\CC^k$-sections
of the bundle $E^*\otimes\vo^{1-q}(X)$, denoted by $\Ga^k_c(X,E^*\otimes\vo^{1-q}(X))$.

To set up an appropriate topology in that space we denote the bundle $E^*\otimes\Vol^{1-q}$ 
by $F$ and define for any $K\comp X$ the space $\Ga^k_{c,K}(X,F):=\{u\in\Ga^k(X,F)\mid
\supp(u)\subseteq K\}$. $\Ga^k_{c,K}(X,F)$ is a Fr\'echet space
and we endow $\Ga^k_c(X,F)$ with the {\em inductive limit topology} $\tau$ with respect to 
the spaces $\Ga^k_{c,K}(X,F)$. $(\Ga^k_{c}(X,F),\tau)$ is isomorphic to the topological
direct sum of the (LF)-spaces $\Ga_c^k(X_\lambda,F)$ where $(X_\la)_\la$ denotes
the family of connected components of $X$ (hence it is an (LF)-space only if $X$ is second countable). 
In particular, $\tau$ is Hausdorff
and complete and we can write $\Ga^k_c(X,F)\,=\,\lim\limits_{\longrightarrow}\Ga^k_{c,K}(X,F)\,.$
This renders the space of compactly supported $\CC^k$-sections of $F$ a {\em
strict inductive limit} (of (F)-spaces), in the sense of Def.\ 2 in Ch.\ 4, Part 1, Sec.\ 3
of \cite{groth}. Most of the typical properties
of strict (LF)-spaces (cf.\ \cite{Schaefer}, 6.4--6.6;
\cite{MK}, 7.1.4) carry over to $\Ga^k_{c}(X,F)$:
It is a complete locally convex space; on each $\Ga^k_{c,K}(X,F)$, $\tau$ induces the
Fr\'echet topology; every bounded subset of $\Ga^k_{c}(X,F)$ is contained (and bounded)
in the Fr\'echet space $\Ga^k_{c,K}(X,F)$ for some $K\comp X$.
 
Finally, the space $\D'\,^{(k)}(X,E\otimes\Vol^q(X))$ of $E$-valued distributions 
of order $k$ and density character $q$ is defined as the topological dual of 
$\Ga_c^k(X,E^*\otimes\Vol^{1-q})$, i.e.,
\begin{equation} \label{distdef}
  \D'\,^{(k)}(X,E\otimes\Vol^q(X)):=[\Ga_c^k(X,E^*\otimes\Vol^{1-q}(X))]'.
\end{equation}
Analogous to the theory on open sets of Euclidean space
the space of smooth regular objects, i.e., $\Ga^\infty(X,E\otimes\Vol^q(X))$
is sequentially dense in $\D'\,^{(k)}(X,E\otimes\Vol^q(X))$.   

We explicitly mention the following special cases of (\ref{distdef}) 
(already anticipated in the discussion above): for 
$E=X\times \C$, $k=\infty$, $q=0$ resp.\ $q=1$ we obtain $\D'(X)$ resp.\ $\D'_d(X)$, 
the space of distributions resp.\ distributional densities on $X$. Similarly, taking
$E$ the tensor bundle $T^r_s(X)$, $k=\infty$ and $q=0$ resp.\ $q=1$ gives the 
spaces $\D'^r_s(X)$ of tensor distributions resp.\  ${\D'_d}{}^r_{s}(X)$ of tensor
distribution densities.

$E$-valued distributions of density character $q$ 
may be written as classical sections 
of $E$ with distributional coefficient ``functions'',
more precisely 
\[
   \D'(X)\otimes_{\CC^\infty(X)}\Ga(X,E\otimes\Vol^q(X))\,
   \,\cong\D'(X,E\otimes\Vol^q(X)).
\]

For $X$ an oriented manifold whose orientation is induced by a fixed nowhere vanishing
$\theta \in \Om^n(X)$,
a rich theory of distributional geometry was introduced by Marsden in
\cite{marsden}. The basic idea underlying his approach is that of continuous extension
of classical operations to spaces of currents: Since $X$ is oriented we may identify
one-densities and smooth $n$-forms and we set
$$
\Om^k(X)' := \D'(X,E\otimes\mbox{Vol}(X))
$$
where $E^*=\Lambda^{n-k}T^*X$.
Using the above identification it follows that $\Om^k(X)'$ is the dual of 
$\Om_c^{n-k}(X)$, the space of compactly supported $n-k$-forms 
(and {\em not} the dual of $\Om^k(X)$ as might be suggested by this notation).
Also, $\D'(X)\cong \Om^0(X)'\cong \D'_d(X)$ and
$\Om^k(X)'$ is precisely the space of odd $k$-currents on $X$ in the
sense of de Rham (\cite{deR}).
Marsden calls elements of $\Om^k(X)'$ {\em generalized $k$-forms} but we prefer
here the term {\em distributional $k$-forms} since the term ``generalized'' will
be reserved for Colombeau objects in this work. Embedding of regular objects into
distributional $k$-forms is effected by the map
\begin{equation}\label{marsdenemb}
\begin{array}{rcl}
j: \Om^k(X) &\to& \Om^k(X)' \\
j(\om)(\tau) &=& \int \om\wedge\tau 
\end{array}
\end{equation}
It then follows that $\Om^k(X)'$ is the weak sequential closure of $j(\Om^k(X))$
(in fact, Marsden {\em defines} $\Om^k(X)'$ as this closure). Let us exemplify
the method of continuously extending classical operations from smooth to distributional
forms by considering the Lie derivative with respect to a smooth vector field $\xi$.
By Stokes' theorem, for $\om\in \Om^k(X)$, $\tau \in \Om_c^{n-k}(X)$ we have
$j(L_\xi\om)(\tau) = - \om(L_\xi\tau)$. Hence setting $L_\xi\om(\tau) := 
- \om(L_\xi\tau)$ for $\om \in \Om^k(X)'$ gives the unique continuous extension 
of $L_\xi$ to $\Om^k(X)'$. By the same strategy, operations like exterior 
differentiation
$d$ and insertion $i_\xi$ can be extended to distributional forms while preserving
classical relations like $L_\xi = i_\xi\circ d + d\circ i_\xi$.
Finally, we note that in this setting, $\D'^r_s(X)$ can be identified with the space
of ${\cal C}^\infty$-multilinear maps $t: {\Om}^1(X)^r\times {\mathfrak X}(X)^s
\to \D'(X)$.

\section{Basic properties, point value characterization} \label{specalgmf}
\blem \label{emlem} Set $\es(X):=({\cal C}^\infty(X))^I$.
The following spaces of nets are equal
\begin{itemize}
\item[(i)]
$\{(u_\eps)_{\eps\in I}\in {\cal E}(X) |\ 
\forall K\subset\subset X,\ \forall P\in{\cal P}(X)\ \exists N\in\N:$\\ 
\hspace*{6.3cm}$\sup_{p\in K}|Pu_\eps(p)|$ 
 $=$ $O(\eps^{-N})\}$
\item[(ii)]
$\{(u_\eps)_{\eps\in I}\in {\cal E}(X) |\
\forall K\subset\subset X,\ \forall k\in\N_0\ \exists N\in\N $\\
\hspace*{1.2cm}$\forall \xi_1,\dots,\xi_k\in {\mathfrak X}(X):\ 
\sup_{p\in K}|L_{\xi_1}\dots L_{\xi_k}\,u_\eps(p)|=O(\eps^{-N})\}$
\item[(iii)] 
$\{(u_\eps)_{\eps\in I}\in {\cal E}(X) | \mbox{ for each chart} \
(V,\psi): $\\ 
\hspace*{6.9cm}$ (u_\eps\circ\psi^{-1})_\eps \in{\cal E}_M(\psi(V))\}$   
\end{itemize}
\elem
\pr Since every iterated Lie derivative is an element of ${\cal P}(X)$ we have
(i) $\subseteq$ (ii). (ii) $\subseteq$ (iii) is immediate from the local form
of $L_{\xi_1}\dots L_{\xi_k}$. Finally,  (iii) $\subseteq$ (i) follows from
Peetre's theorem (see e.g., \cite{Kahn}, Th.\ 6.2). \ep

We denote by ${\cal E}_M(X)$ the set defined above and call it the space of
{\em moderate} nets on $X$. Definition (i) was suggested in \cite{RD},
(iii) is from \cite{AB}. (ii) is mentioned explicitly since the operation of
taking Lie derivatives plays a central role in the theory (in the full version
of the construction, a canonical embedding of $\D'$ commuting with Lie derivatives
has been given in \cite{vim}). Replacing $\exists N$ by $\forall m$, and 
$\eps^{-N}$ by $\eps^m$ in (i) and (ii) as well as ${\cal E}_M(\psi(V))$ by
${\cal N}(\psi(V))$ in (iii) we obtain equivalent definitions of the space
${\cal N}(X)$ of {\em negligible} nets on $X$. Applying \cite{found},
Th.\ 13.1
locally, we arrive at the following characterization of ${\cal N}(X)$ as a 
subspace of ${\cal E}_M(X)$:
\begin{equation} \label{nxnew}
{\cal N}(X) = \{(u_\eps)_\eps\in {\cal E}_M(X) \mid \forall K\comp X \, \forall
m\in \N \, \sup_{x\in K} |u_\eps| = O(\eps^m)\}
\end{equation}
Thus for elements of ${\cal E}_M(X)$ to belong to ${\cal N}(X)$ it suffices 
to require the ${\cal N}$-estimates to hold for the function itself, without
taking into account any derivatives.
The {\em Colombeau algebra of generalized functions on the manifold $X$}
is defined as the quotient 
\[ 
        \gs(X)\,:=\,\esm(X)\,/\,\ns(X)\,.
\]

Again, elements in $\gs(X)$
are denoted by capital letters, i.e., $U=\cl[(u_\eps)_\eps]=
(u_\eps)_\eps+\ns(X)$. Analogous to the case of open sets in 
Euclidean space, $\esm(X)$ is a differential
algebra (w.r.t.\ Lie derivatives) 
with componentwise operations and $\ns(X)$ is a differential ideal in it. 
Moreover, ${\cal E}_M(X)$ and ${\cal N}(X)$ are invariant under 
the action of any $P\in {\cal P}(X)$. Thus we obtain
\bp
Let $U\in\gs(X)$ and $P\in {\cal P}(X)$. Then 
\[ P U\,:=\,\cl[(Pu_\eps)_\eps]\]
is a well-defined element of $\G(X)$
\eprop
This applies, in particular, to the Lie derivative $L_\xi U$ of $U$ with respect
to a smooth vector field $\xi \in {\mathfrak X}(X)$.
It follows that $\gs(X)$ is a differential $\K$-algebra w.r.t. Lie derivatives. 

It is now immediate that a generalized function $U$ on $X$ allows for the following 
local description via the assignment $\gs(X)\ni U\mapsto (U_\al)
_{\al\in A}$ with $U_\al:=U\circ\psi_\al^{-1}\in\gs(\psi_\al(V_\al))$
(with $\{(V_\al,\psi_\al)\mid \al\in A\}$ an atlas of $X$).
We call $U_\al$ the {\em local expression} of $U$ with respect to the 
chart $(V_\al,\psi_\al)$. Thus we have
\bp\label{glocal} $\gs(X)$ can be identified with the set of all families $(U_\al)_\al$ 
of generalized functions $U_\al\in\gs(\psi_\al(V_\al))$ satisfying the following transformation law
\[
        U_\al|_{\psi_\al(V_\al\cap V_\beta)}
        \,=\,
        U_\beta|_{\psi_\beta( V_\al\cap V_\beta)}\circ\psi_\beta\circ\psi_\al^{-1}
\]
for all $\al,\beta\in A$ with $V_\al\cap V_\beta\not=\emptyset$. 
\ep
\eprop

It follows that $\gs(\_)$ is a fine sheaf of $\K$-algebras on $X$. In fact,
in \cite{RD}, $\G$ is {\em defined} directly as a quotient sheaf of the sheaves
of moderate modulo negligible sections. 

An important feature distinguishing Colombeau generalized functions on open
subsets $\Om$ of $\R^n$ from spaces of distributions is the availability of
a point value characterization  of elements of $\G(\Om)$ (\cite{point}). This
characterization allows a direct generalization of results from classical
analysis to Colombeau algebras thereby enabling a consistent treatment of a
variety of geometric and analytic problems (see e.g., \cite{HO},
\cite{symm}).
Our aim in the remainder of this section is to derive a point value characterization 
of Colombeau generalized functions also in the global context.

To begin with we shortly recall the basic notions from \cite{point}. 
Clearly a generalized function is not characterized by its values
on all classical points: on $\R$, take $F=\iota(x)\iota(\de)$; then
$F\not=0$ but $F(x)=0$ in ${\cal R}$ $\forall x\in\R$. The basic idea is therefore to 
introduce an analogue of ``nonstandard numbers'' into the theory which are flexible enough to
capture all the relevant information contained in a generalized function. 
Let $\Om\subseteq \R^n$ open. We define the set of compactly supported
sequences of points on $\Om$ by $\Om_c := \{(x_\eps)_\eps\in \Om^I \mid \exists
K\subset\subset\Om \mbox{ such that } x_\eps\in K\ \forall\eps\mbox{ small}\}$. 
Next we introduce the following equivalence relation: two elements 
$(x_\eps)_\eps$, $(y_\eps)_\eps\in\Om_c$ are called equivalent 
($(x_\eps)_\eps\sim (y_\eps)_\eps$) if $|x_\eps -y_\eps| = O(\eps^m)$ for each $m>0$. 
Finally we define  the {\em set of compactly
supported generalized points} as the quotient $\widetilde\Om_c:=\Om_c/\!\!\sim\,$.
Then for any $U\in\G(\Om)$ and $\tilde x \in \widetilde\Om_c$, the generalized
point value $U(\tilde x):=\cl[(u_\eps(x_\eps))_\eps]$
is a well-defined generalized number (\cite{point}, Prop.\ 2.3). Moreover,
generalized functions on $\Om$ are characterized by their generalized point values
in the sense that $U=0$ iff
$U(\tilde x) = 0$ for each $\tilde x$ in $\widetilde\Om_c$ (\cite{point},
Th.\ 2.4).

In order to transfer these notions to the manifold-setting we will make use
of
an auxiliary Riemannian metric $h$ on $X$. Of course we will then have to
show that
the constructions to follow are in fact independent of the chosen $h$.

We call a net $(p_\eps)_\eps \in X^I$ {\em compactly supported} if
there exist
$K\comp X$ and $\eta>0$ such that $p_\eps\in K$ for $\eps<\eta$. Denoting by
$d_h$ the Riemannian distance induced by $h$ on $X$, two nets $(p_\eps)_\eps$,
$(q_\eps)_\eps$
are called equivalent ($(p_\eps)_\eps\sim (q_\eps)_\eps$) if
$d_h(p_\eps,q_\eps)
=O(\eps^m)$ for each $m>0$. The equivalence classes with respect to this
relation
are called {\em compactly supported generalized points} on $X$. The set of
compactly supported generalized points on $X$ will be denoted by $\widetilde
X_c$.

The fact that $\widetilde X_c$ does not depend on the auxiliary metric $h$ 
follows immediately from the following lemma:

\blem\label{riemlem} Let $h_i$ be Riemannian metrics inducing the Riemannian
distances $d_i$ on $X$
($i=1,2$).
Then for $K, K' \comp X$ there exists $C>0$ such that $d_2(p,q) \le C
d_1(p,q)$
for all $p\in K$, $q\in K'$.
\elem
\pr Assume to the contrary that there exist sequences $p_m$ in $K$ and $q_m$
in
$K'$ such that $d_2(p_m,q_m) > m d_1(p_m,q_m)$. By choosing suitable
subsequences
we may additionally suppose that both $p_m$ and $q_m$ converge to some $p$.
Let $V$ be a relatively compact neighborhood of $p$. Then denoting by $B_r^i(q)$
the $d_i$-ball of radius $r$ around $q$ it follows that there exist $r_0>0$ and
$\alpha>0$ such that $B_r^1(q) \subseteq  B_{\alpha r}^2(q)$ for all $q\in V$ and
all $r<r_0$ (cf.\ e.g., \cite{vim}, Lemma 3.4). But then for $m>\alpha$ sufficiently
large we arrive at the contradiction $d_2(p_m,q_m) \le \alpha d_1(p_m,q_m)$.\ep

\blem \label{simlocal}
Suppose that $(p_\eps)_\eps$, $(q_\eps)_\eps \in X^I$ are  compactly
supported in some $W_\al$ which is open, geodesically convex 
with respect to a Riemannian metric $h$ on $X$
and
satisfies $\overline{W_\al}\comp V_\al$  for some chart $(V_\al,\psi_\al)$.
Then 
$$
d_h(p_\eps,q_\eps)=O(\eps^m)\ \forall m>0\Leftrightarrow 
|\psi_\al(p_\eps) - \psi_\al(q_\eps)| = O(\eps^m) \ \forall \, m>0\,.
$$
\elem
\pr $(\Rightarrow)$  
Let $\gamma_\eps : [\al_\eps,\beta_\eps] \to W_\al$ be the unique geodesic in $W_\al$
joining $p_\eps$ and $q_\eps$. Then 
$$
d_h(p_\eps,q_\eps) = \int_{\al_\eps}^{\beta_\eps} \|\gamma_\eps'(s)\|_h \,ds
= O(\eps^m) \ \forall \, m>0\,.
$$
Since $W_\al$ is relatively compact there exists $C>0$ such that 
$|\xi|\le C\,\|T_{\psi_\al(p)}\psi_\al^{-1}\xi\|_h$ for all 
$p\in W_\al$ and all $\xi\in\R^n$.
Thus 
$$
|\psi_\al(p_\eps) - \psi_\al(q_\eps)| \le 
\int_{\al_\eps}^{\beta_\eps} |(\psi_\al\circ\gamma_\eps)'(s)| \,ds
\le C \int_{\al_\eps}^{\beta_\eps} \|\gamma_\eps'(s)\|_h \,ds = O(\eps^m)
$$

$(\Leftarrow)$ Let $K\comp W_\al$ such that $p_\eps, q_\eps \in K$ for $\eps$
small. Using a cut-off function supported in $\psi_\al(V_\al)$  
and equal to $1$ in a neighborhood $W'$ with $\overline{W'}\comp \psi_\al(W_\al)$ of 
$\psi_\al(K)$ we may extend the pullback under $\psi_\al$ 
of the Euclidean metric on $\psi_\al(V_\al)$ to a Riemannian metric $g$ on $X$.
There exists $\eps_0>0$ such that for each $\eps<\eps_0$ the whole line 
connecting $\psi_\al(p_\eps)$ with $\psi_\al(q_\eps)$ is contained in $W'$. Hence
$d_g(p_\eps,q_\eps) \le d_{g|_{\psi_\al^{-1}(W')}}(p_\eps,q_\eps) = 
|\psi_\al(p_\eps)-\psi_\al(q_\eps)| = O(\eps^m)$,
so the claim follows from Lemma \ref{riemlem}.\ep



\bp \label{pointwelldef}
Let $U\in \G(X)$ and $\tilde p \in \widetilde X_c$. Then $$ U(\tilde p):=
\cl[(u_\eps(p_\eps))_\eps]$$ is a well-defined element of ${\cal K}$. 
\eprop
\pr Since $(p_\eps)_\eps$ is compactly supported it is clear that $(u_\eps(p_\eps))$
is moderate resp.\ negligible if $(u_\eps)_\eps$ is.  
Suppose now that $(p_\eps)_\eps \sim (q_\eps)_\eps$ and choose
$K\comp X$ such that $p_\eps, q_\eps \in K$ for $\eps$ small. We have to show that
$(u_\eps(p_\eps) - u_\eps(q_\eps))_\eps \in {\cal N}$.
To this end we choose some
auxiliary Riemannian metric $h$ and cover $K$ by finitely many $W_{\al_i}$
with $\overline{W_{\al_i}}\comp V_{\al_i}$ as in Lemma \ref{simlocal}.
$K$ can be written as the union of 
compact sets $K_i \comp W_{\al_i}$. 
Then for each $\eps$ sufficiently small there exists $i_\eps$ such that 
the line connecting $\psi_{\al_{i_\eps}}(p_\eps)$ with $\psi_{\al_{i_\eps}}(q_\eps)$ 
is contained in $\psi_{\al_{i_\eps}}(W_{\al_{i_\eps}})$.
Thus the claim follows from Lemma \ref{emlem} 
and Lemma \ref{simlocal} 
by applying the mean value theorem as in \cite{point}, Prop.\ 2.3.\ep 

\btheorem Let $U\in \G(X)$. Then $U = 0$ in $\G(X)$ iff $U(\tilde p) = 0$ in ${\cal K}$
for all $\tilde p \in \widetilde X_c$.
\etheorem
\pr Necessity is immediate from Proposition \ref{pointwelldef}. 
Conversely, fix some Riemannian
metric $h$ and cover $X$ by geodesically convex sets $W_\al$ 
with $\overline{W_\al}\comp V_\al$ for charts 
$(V_\al,\psi_\al)$. Let $\tilde x \in \psi_\al(W_\al)^\sim_c$. Then by Lemma 
\ref{simlocal}
$\tilde p := \cl[(\psi_\al^{-1}(x_\eps))_\eps]$ is a well-defined element of 
$\widetilde X_c$. By assumption $u_\eps(p_\eps) = 
u_\eps\circ\psi_\al^{-1}(x_\eps)$ is a negligible net in $\K$. Thus by 
\cite{point}, Th.\ 2.4, $U\circ\psi_\al^{-1} = 0$ in $\G(\psi_\al(W_\al))$ 
for all $\al$, so $U=0$ by Proposition \ref{glocal}.\ep


\section{Compatibility with distributional geometry, embeddings, and association}
\label{assoz1}
As in \cite{RD} we call $U\in \G(X)$ {\em associated} to $0$, $U\approx 0$, if
$\int_X u_\eps \mu \to 0$ ($\eps\to 0$) for all compactly supported one densities
$\mu \in \Gamma_c^\infty(X,\mbox{Vol}(X))$ and one (hence every) representative
$(u_\eps)_\eps$ of $U$. Clearly, $\approx$ induces an equivalence relation on
$\G(X)$ giving rise to a linear quotient space. 
If $\int_X u_\eps \mu \to w(\mu)$ for some $w\in\D'(X)$ then $w$ is
called the {\em distributional shadow} (or {\em macroscopic aspect}) of $U$
and we write $U\approx w$.
In terms of the local description established in Proposition \ref{glocal} we have
\begin{equation}
\label{wwi1}
       U\approx 0  
       \,\Leftrightarrow\, U_\al\approx 0 \mbox{ in }\gs(\psi_\al(V_\al))\quad\forall\al
\end{equation}
From this it follows that $U_1\approx U_2$ implies $PU_1 \approx PU_2$ for each
$P\in {\cal P}(X)$.

By \cite{hoer1}, 6.3.4, any $w\in \D'(X)$ can be identified with a family $(w_\al)_{\al\in A}$,
where $w_\al\in \D'(\psi_\al(V_\al))$ satisfies the transformation law
$$
w_\beta = (\psi_\al\circ\psi_\beta^{-1})^*(w_\al)\,.
$$
Here $f^*w$ denotes the pullback of a distribution $w$ under the diffeomorphism $f$.
In particular, $w_\al = (\psi_\al^{-1})^*(w|_{V_\al})$. Again a straightforward
calculation gives
\begin{equation}
       U\approx w  
       \,\Leftrightarrow\, U_\al\approx w_\al \mbox{ in }\gs(\psi_\al(V_\al))\quad\forall\al
\end{equation}
Association relations will be our main tool in establishing compatibility with linear
distributional geometry later on. Before we proceed with this analysis, however, let us
address the problem of embedding ${\cal C}^\infty(X)$ and $\D'(X)$ into $\G(X)$.
As in the case of open subsets of $\R^n$, ${\cal C}^\infty(X)$ is embedded into
$\G(X)$ via the ``constant'' embedding $\sigma:\CC^\infty(X)\hookrightarrow\gs(X)$, 
$f \mapsto \cl[(f)_\eps]$.

Turning now to the interrelation between $\D'(X)$ and $\gs(X)$ let us first clarify
what we can expect at all from such an embedding. The method of choice for open subsets
of $\R^n$, i.e., convolution with a mollifier $\rho$ as in Section \ref{intro} is manifestly 
not diffeomorphism invariant, as is demonstrated by the following simple

\bex
Consider the diffeomorphism $\mu(x)=2x$ on $\R$ and set $w=\de \in \D'(\R)$.
Then $\mu^*\de=\frac{1}{2}\de$ and we have
\beas
        ((\iota\circ\mu^*)\,\de)_\eps&=&\iota(\frac{1}{2}\,\de)_\eps\,
        =\,\frac{1}{2}\,\rho_\eps\nn\\
        ((\mu^*\circ\iota)\,\de)_\eps&=&\mu^*\rho_\eps\,=\,\rho_\eps(2\,.)\,.
\eeas
From this we see that $((\iota\circ\mu^*-\mu^*\circ\iota)\de)_\eps=
\frac{1}{2}\rho_\eps(x)-\rho_\eps(2x)$ is not in the ideal $\ns(\R)$. 
However, it is evident that 
$(\iota\circ\mu^*-\mu^*\circ\iota)\de\approx 0$. In fact diffeomorphism invariance
does hold on the level of association (cf.\cite{AB}, Th.\ 9.1.2). 
\eex 
Finally, as was shown in \cite{RD}, Remark 3, there can be no embedding of $\D'(X)$ into
$\G(X)$ that commutes with differentiation in all local coordinates. The fact that a
canonical embedding commuting with Lie derivatives was constructed in \cite{vim} for
the full Colombeau algebra rests heavily on the dependence of representatives
on an additional parameter $\phi \in \D(X)$ (and on the ensuing modified definition
of Lie derivatives of such representatives). Therefore we cannot expect an embedding
providing this property in the setting of the special Colombeau algebra on manifolds.

On the positive side, the existence of injective sheaf morphisms $\iota:
\D'\hookrightarrow\gs$ coinciding with $\sigma$ on $\CC^\infty$ and 
satisfying $\iota(w) \approx w$ for each $w \in\D'(X)$ 
has been proved by de Roever and Damsma~\cite{RD} using de 
Rham-regularizations (cf.\ \cite{deR}, \S 15). In view of the above restrictions these
properties of the embedding seem optimal (unless one is willing to furnish
$X$ with additional structure).

In the following construction\footnote{suggested by M. Oberguggenberger} 
we give an embedding which, while also providing a sheaf 
morphism possessing these optimal properties, is considerably simpler than the 
construction in \cite{RD}, Th.\ 1. 

\btheorem \label{moemb}
Let ${\cal A}=(\psi_\alpha,V_\alpha)_\alpha$ be an atlas of $X$ and let
$\{\chi_j:j\in \N\}$ a smooth partition of unity subordinate to $(V_\alpha)
_\alpha$. Let $\mbox{supp}(\chi_j)\subseteq V_{\alpha_j}$ for $j\in \N$
and choose for every $j\in \N$ some $\zeta_j\in {\cal D}(V_{\alpha_j})$ such
that $\zeta_j\equiv 1$ on $\mbox{supp}(\chi_j)$. Fix some mollifier
$\rho\in {\cal S}(\R^n)$ with unit integral and $\int \rho(x) x^\al\,dx=0$
for all $|\al| \ge 1$.
The map
\beast
& \iota_{\cal A}: {\cal D}'(X) \to {\cal G}(X) &\\
& u \to \cl [
(\sum\limits_{j=1}^\infty
\zeta_j \cdot (((\chi_j \circ \psi_{\alpha_j}^{-1})u_{\alpha_j})
\ast \rho_\eps)\circ \psi_{\alpha_j})_\eps]     &
\eeast
is a linear embedding that coincides with $\sigma$ on ${\cal C}^\infty(X)$.
Moreover, for each $u\in \D'(X)$ we have $\iota_{\cal A}(u) \approx u$ 
and $\supp(u) = \supp(\iota_{\cal A}(u))$.
\etheorem
\pr In the proof we will for the sake of brevity replace $\alpha_j$
by $j$ and set $\widetilde{V}_\alpha = \psi_\alpha(V_\alpha)$. It is obvious
that 
\[
u_\eps := \sum\limits_{j=1}^\infty
\zeta_j \cdot (((\chi_j \circ \psi_j^{-1})u_j)\ast 
\rho_\eps)\circ \psi_j
\]
is a smooth function on $X$. Our first task will therefore consist in 
verifying the ${\cal E}_M$-bounds for $(u_\eps)_\eps$. This means that we
have to estimate $u_\eps\circ \psi_\alpha^{-1}$ for arbitrary $\alpha\in A$.
Let $K\subset\subset \widetilde{V}_\alpha$. Then $L_j = 
\psi_j(\mbox{supp}(\zeta_j)\cap \psi_\alpha^{-1}(K))$ is a compact subset
of $\widetilde{V}_j$. The fact that the 
${\cal E}_M(\widetilde{V}_j)$-function 
$(((\chi_j \circ \psi_j^{-1})u_j)\ast \rho_\eps)_\eps$ 
satisfies the necessary bounds on $L_j$ shows that 
$(u_\eps\circ \psi_\alpha^{-1})_\eps \in {\cal E}_M(\widetilde{V}_\alpha)$. 

To prove injectivity of $\iota_{\cal A}$, we suppose that 
$(u_\eps\circ \psi_\alpha^{-1})_\eps \in {\cal N}(\widetilde{V}_\alpha)$
for all $\alpha\in A$. We have to show that $u_\alpha = 0$ in
${\cal D}'(\widetilde{V}_\alpha)$ for all $\alpha$. Fix some $\alpha\in A$
and let $\varphi\in {\cal D}(\widetilde{V}_\alpha)$. The
term $\langle u_\eps\circ \psi_\alpha^{-1},\varphi\rangle$ is a finite sum
of expressions of the form
\beast
& \int\limits_{\widetilde{V}_\alpha} \zeta_j\circ \psi_\alpha^{-1}(x)
(((\chi_j \circ \psi_j^{-1})u_j)\ast \rho_\eps)(\psi_j\circ\psi_\alpha^{-1})
(x)\varphi(x)\,dx& \\
& =\int\limits_{\psi_\alpha(V_j\cap V_\alpha)} 
\zeta_j\circ \psi_\alpha^{-1}(x)
(((\chi_j \circ \psi_j^{-1})u_j)\ast \rho_\eps)(\psi_j\circ\psi_\alpha^{-1})
(x)\varphi(x)\,dx&\\
& =\int\limits_{\psi_j(V_j\cap V_\alpha)} 
\zeta_j\circ \psi_j^{-1}(y)
(((\chi_j \circ \psi_j^{-1})u_j)\ast \rho_\eps)(y)\varphi
\circ \psi_\alpha\circ\psi_j^{-1}(y) &\\
& |\mbox{det}(D(\psi_\alpha\circ\psi_j^{-1}))(y)|\,dy.&
\eeast
For $\eps\to 0$, this converges to
\beast
\lefteqn{\langle \zeta_j\circ \psi_j^{-1}\cdot (\chi_j\circ\psi_j^{-1})u_j,
\varphi\circ \psi_\alpha\circ\psi_j^{-1}
|\mbox{det}(D(\psi_\alpha\circ\psi_j^{-1}))|\rangle
}  \\
&& =\langle \zeta_j\circ\psi_\alpha^{-1}\circ\psi_\alpha\circ 
\psi_j^{-1}\cdot (\chi_j\circ\psi_\alpha^{-1}\circ\psi_\alpha\circ
\psi_j^{-1})u_j,\\
&& \hphantom{mmmmmmmmmmmm}
\varphi\circ \psi_\alpha\circ\psi_j^{-1}
|\mbox{det}(D(\psi_\alpha\circ\psi_j^{-1}))|\rangle  \\
&& 
=\langle (\zeta_j\circ \psi_\alpha^{-1})(\chi_j\circ \psi_\alpha^{-1})
(\psi_j\circ\psi_\alpha^{-1})^\ast u_j, \varphi \rangle \\
&& =\langle (\zeta_j\circ \psi_\alpha^{-1})(\chi_j\circ \psi_\alpha^{-1})
u_\alpha, \varphi \rangle. 
\eeast
Therefore, for $\eps\to 0$ we have
\beast
& \langle u_\eps\circ \psi_\alpha^{-1},\varphi\rangle
\to \sum\limits_{j=1}^\infty \langle 
(\zeta_j\circ \psi_\alpha^{-1})(\chi_j\circ \psi_\alpha^{-1})
u_\alpha, \varphi \rangle &\\
& =\sum\limits_{j=1}^\infty \langle 
(\chi_j\circ \psi_\alpha^{-1})
u_\alpha, \varphi \rangle = \langle u_\alpha, \varphi \rangle. &
\eeast
On the other hand, since $(u_\eps)_\eps\in {\cal N}(X)$, the above
expression converges to $0$, which establishes the injectivity of 
$\iota_{\cal A}$. Also, the above calculation shows that
$\iota_{\cal A}(u) \approx u$ for each $u\in \D'(X)$.

Let $f\in {\cal C}^\infty(X)$. We claim that 
$U:=\iota_{\cal A}(f)=\sigma(f)$. Considered as an element of ${\cal D}'(X)$,
$f$ is identified with $((f\circ \psi_\alpha^{-1})_\alpha)_\alpha$, so
\beast
& u_\eps = \sum\limits_{j=1}^\infty
\zeta_j \cdot (((\chi_j f)\circ \psi_j^{-1})\ast 
\rho_\eps)\circ \psi_j.&
\eeast
We have to show that $((u_\eps-f)\circ\psi_\alpha^{-1})_\eps\in {\cal N}
(\widetilde{V}_\alpha)$ for all $\alpha\in A$. Now
\beast
& f(x)=\sum\limits_{j=1}^\infty \zeta_j(x)(\chi_j\cdot f)(x) =
\sum\limits_{j=1}^\infty \zeta_j(x)((\chi_j\cdot f)\circ \psi_j^{-1})
(\psi_j(x)), &
\eeast
so
\beast
&&(u_\eps - f)\circ \psi_\alpha^{-1}\\
&&\quad=\sum\limits_{j=1}^\infty \zeta_j\circ\psi_\alpha^{-1}\underbrace{
[(((\chi_j\cdot f)\circ \psi_j^{-1})\ast \rho_\eps) - 
(\chi_j\cdot f)\circ \psi_j^{-1}]}_{(\ast)} 
\circ \psi_j\circ \psi_\alpha^{-1}.
\eeast
It therefore suffices to notice that each of the terms $(\ast)$ is in
${\cal N}(\widetilde{V}_j)$. But this follows by Taylor expansion as
in the corresponding proof for open subsets of $\R^n$.
Finally, preservation of supports is also deduced exactly as in the
local case (cf.\ e.g., \cite{kunzDISS}, 1.2.8).
\ep

It immediately follows that $\iota_{\cal A}$ is
a local operator, i.e., it indeed induces a sheaf morphism with the 
above properties. Nevertheless, just as the corresponding construction
in \cite{RD} $\iota_{\cal A}$ is {\em non-geometric} in an essential way,
i.e., it depends on the chosen atlas as well as on the functions $\zeta_j$,
$\chi_j$, etc. For practical purposes however, this drawback is often compensated
by the availability of regularization procedures adapted to the specific problem
at hand that can be used to model the singularities directly in $\G(X)$ 
without the use of a distinguished embedding. The connection to the distributional
picture is then effected by means of association procedures (cf.\ e.g., 
\cite{ultra}, \cite{penrose}) whose basic properties we now continue to study.

To this end let us first discuss consistency properties with respect to classical 
products (in the sense of association). 
In the absence of a distinguished embedding $\iota$ we have  to be slightly 
more cautious than in the case of $\R^n$.
For example the following (naive) generalization of the  statement that
the product $\CC^\infty\times\D' \to \D'$ is respected by association
(more precisely $\iota(f)\iota(u)\approx\iota(fu)$ for all
$f\in\CC^\infty(\Om),\ u\in\D'(\Om)$): ``$U,V\in\gs(X)$, 
$U\approx f\in\CC^\infty$ and $V\approx w\in\D'(X)$ $\Rightarrow$
$UV\approx fw$'' is wrong in general. 
To see this take $\rho\in\D(\R)$ with $\int\rho=1$. 
Then $\cl[(\rho(\frac{x}{\eps}))_\eps]\approx 0$ and clearly $\cl[(\frac{1}{\eps})
\rho(\frac{x}{\eps})_\eps]\approx\de$
but $\rho(\frac{x}{\eps})\,(\frac{1}{\eps})\rho(\frac{x}{\eps})\to
\de \int\rho^2$ in $\D'$.  
The reason for the validity of the corresponding $\R^n$-statement 
ultimately is that $f*\rho_\eps\to f$
uniformly on compact sets already for a continuous function $f$, whereas 
$\rho(x/\eps)\to 0$ only weakly.
Therefore we introduce the following stronger equivalence relations on $\gs(X)$.

\bd 
Let $U\in\gs(X)$.\begin{itemize}
\item[(i)] $U$ is called {\em $\CC^k$-associated} to $0$
$(0\leq k\leq\infty)$, $U\approx_k 0$, if for all $l\leq k$, 
all $\xi_1,\dots,\xi_l
\in {\mathfrak X}(X)$ and one (hence any) representative $(u_\eps)_\eps$
\[ L_{\xi_1}\dots L_{\xi_l}\,u_\eps\,\to\,0\, \mbox{ uniformly on compact sets.}
\] 
\item[(ii)] We say that $U$ admits $f$ as {\em $\CC^k$-associated
function}, 
$U\approx_k f$, if for all $l\leq k$,
all $\xi_1,\dots,\xi_l \in {\mathfrak X}(X)$ and one (hence any) representative
\[ L_{\xi_1}\dots L_{\xi_l}\,(u_\eps-f)\,\to\,0 \mbox{ uniformly on compact sets.}
\] 
\end{itemize}
\edefi

Clearly if $U$ is $\CC^k$-associated to $f$ then $f\in\CC^k(X)$. Moreover, if
$U$ admits for a $\CC^k$-associated function at all the latter is unique. 
Note also that the above notion of convergence may equivalently be expressed by 
saying that all $(u_\al\,_\eps)_\eps$ converge uniformly in all derivatives of order 
less or equal $k$ (resp. in all derivatives if $k=\infty$)  
on compact sets. We are now prepared to state the following
\bp\label{consist} Let $U,V\in\gs(X)$.
\begin{itemize}
\item[(i)] If $V\approx w\in\D'(X)$, $f\in\CC^\infty(X)$, and either (a) $U=\sigma(f)$ or
(b) $U\approx_\infty f$, then $UV\approx fw$.
\item[(ii)] If $U\approx_k f$ and $V\approx_k g$ then $UV\approx_k fg$ ($f,g\in
\CC^k(X)$).
\end{itemize}
\eprop
\pr (i)(a) is clear since $\int f v_\eps\mu=v_\eps(f\mu)\to w(f\mu)$
for all compactly supported one-densities $\mu$.
To prove (i)(b) we use the fact that multiplication: $\CC^\infty\times\D'\to\D'$ as 
a bilinear separately continuous map is jointly sequentially continuous 
since both factors are barrelled (\cite{koethe}, \S42.2(3) and \S40.1). 
(ii) follows from elementary analysis.
\ep

Proposition \ref{consist} (i)(a) is the reconciliation of the respective $\CC^\infty$-module
structures of $\D'$ and $\gs$ on the level of association.
Next we introduce the notion of integration of generalized functions.
\bd \label{densityint}
Let $U\in\gs(X)$ and $\mu\in\Ga^\infty(X,\Vol(X))$. Then we define the {\em integral
of $U$ with respect to $\mu$} over $M\subset\subset X$ by
\[
        \int\limits_M U\mu\,=\,\cl[(\int\limits_M u_\eps\mu)_\eps]\,.
\]
\edefi

For $U\mu$ compactly supported we set $\int_X U\mu := \int_K U\mu$ where $K$
is any compact set containing $\supp(U\mu)$ in its interior. It is easily seen
that this definition is independent of the chosen $K$. Also, we have 
$\int_\R \delta(x)\,dx = 1$. We close this section 
by showing that the Lie derivative respects associated distributions.

\bp Let $X$ be orientable and $U\approx w$. Then $L_\xi U\approx L_\xi w$.
\eprop
Orientability is supposed in order to be able to identify one-densities with
$n$-forms, where a Lie derivative is defined.
Moreover, Stokes' theorem is used in the following\\
\pr Let $\nu\in\Om^n_c(X)$ then
\[ 
        \int (L_\xi u_\eps)\nu=-\int u_\eps (L_\xi\nu)\to
        -w(L_\xi\nu)=L_\xi w(\nu)
\]
\ep
\section{Generalized sections of vector bundles}\label{gsvb}

For a section $s\in\Ga(X,E)$ we call $s_\al^i:=\Psi_\al^i\circ s\circ\psi_\al^{-1}$
its $i$-th component with respect to the vector bundle chart $(V_\al,\Psi_\al)$
($i=1,\dots,{n'}$, where ${n'}$ is the dimension of the fibers).

\bd Let $E\to X$ be a vector bundle, and again $I=(0,1]$.
\beas
        \Ga_\es(X,E)&:=& (\Ga(X,E))^I\\
        \Ga_{\esm}(X,E)&:=& \{ (s_\eps)_{\eps\in I}\in \Ga_\es(X,E) : 
                \ \forall\al, \forall i=1,\dots,{n'}:\\
                 &&\hphantom{mmmmm} (s^i_\al\,_\eps)_\eps:=
                (\Psi^i_\al\circ s_\eps\circ\psi_\al^{-1})_\eps
                \in\esm(\psi_\al(V_\al))\}\\ 
        \Ga_{\ns}(X,E)&:=& \{ (s_\eps)_{\eps\in I}\in \Ga_\es(X,E) : 
                \ \forall\al, \forall i=1,\dots,{n'}:\\
                &&\hspace{5.2cm} (s^i_\al\,_\eps)_\eps\in\ns(\psi_\al(V_\al))\}\nn
\eeas
\edefi
First note that although the composition $f\circ U$ of a generalized function $U$ with a 
smooth function $f$ generally need not be moderate the notions of moderateness and 
negligibility as defined above are  preserved under the change of bundle charts due 
to the (fiberwise) linearity of the transition functions. 
In particular, these notions do not depend on the chosen atlas. In fact, using
Peetre's theorem we obtain the following global description of moderate resp.\
negligible sections:
\beas
        \Ga_{\esm}(X,E)&=& \{ (s_\eps)_{\eps\in I}\in \Ga_{\es}(X,E) : 
                \ \forall P\in {\cal P}(X,E)\,\\
                 && \hspace{1cm}
                 \forall K\comp X \, \exists N\in \N:\ \sup_{p\in K}\|Pu_\eps(p)\| = O(\eps^{-N})\}\\ [.5em]
        \Ga_{\ns}(X,E)&=& \{ (s_\eps)_{\eps\in I}\in \Ga_{\es}(X,E) : 
                \ \forall P\in {\cal P}(X,E)\, \\
                 &&\hspace{1.3cm}
                 \forall K\comp X \, \forall m\in \N:\ \sup_{p\in K}\|Pu_\eps(p)\| = O(\eps^{m})\}\\ 
\eeas
Here $\|\, \|$ denotes the norm induced on the fibers of $E$ by any Riemannian
metric. Similar to (\ref{nxnew}), \cite{found}, Th.\ 13.1 yields a characterization 
of $\Ga_{\ns}(X,E)$ as a subspace of 
$\Ga_{{\cal E}_M}(X,E)$ that imposes the above growth restrictions
on representatives only with respect to differential operators of order $0$.
In order to define generalized sections of the bundle $E\to X$ we need the
following 
\bp With operations defined componentwise (i.e., for each $\eps$), 
$\Ga_{\esm}(X,E)$ is a $\gs(X)$-module with $\Ga_{\ns}(X,E)$ a submodule in it.
\eprop
\pr We need to establish the following statements
(a)\,$(u_\eps)_\eps\in\esm(X)$, $(s_\eps)_\eps\in\Ga_{\esm}(X,E)\Rightarrow (u_\eps s_\eps)_\eps\in
\Ga_{\esm}(X,E)$, (b) $(u_\eps)_\eps\in\ns(X)$, $(s_\eps)_\eps\in\Ga_{\esm}(X,E)$ $\Rightarrow$ 
$(u_\eps s_\eps)_\eps\in\Ga_{\ns}(X,E)$ and (c) $(u_\eps)_\eps\in\esm(X)$, $(s_\eps)_\eps\in\Ga_{\ns}(X,E)
\Rightarrow (u_\eps s_\eps)_\eps\in\Ga_{\ns}(X,E)$, which easily follow from the local 
description in Proposition \ref{glocal} and the definitions above.
\ep

Now we are in the position to define.
\bd The $\gs(X)$-module of generalized sections of $E\to X$ is defined as the 
quotient
\[ 
        \Ga_{\ns}(X,E)\,:=\,\Ga_{\esm}(X,E)\,/\,\Ga_{\ns}(X,E)\,.
\]
\edefi

As usual we denote generalized objects by capital letters, e.g.,
$S=\cl[(s_\eps)_\eps]$. 
By the very definition of $\Ga_{\gs}(X,E)$ we may describe a generalized section $S$ by a 
family $(S_\al)_\al=((S^i_\al)_\al)_{i=1}^{n'}$, 
where $S_\al$ is called the 
{\em local expression}
of $S$. Its {\em components} 
$S^i_\al:=\Psi^i_\al\circ S\circ\psi_\al^{-1}
\in\gs(\psi_\al(V_\al))$ ($i=1,\dots, {n'}$) satisfy
\beq\label{gloc}
        S^i_\al(x)\,
=\,(\bfs{\psi}_{\al\beta})^i_j(\psi_\beta\circ\psi^{-1}_\al(x))\,S^j_\beta 
(\psi_\beta\circ\psi^{-1}_\al(x))
\eeq
for all $x\in \psi_\al(V_\al\cap V_\beta)$, where $\bfs{\psi}_{\al\beta}$ denotes
the transition functions of the bundle. Hence formally generalized sections 
of $E\to X$ are locally simply given by ``ordinary'' sections with generalized 
``coefficients.'' We shall see shortly that this property in fact also holds
globally (cf.\ Theorem \ref{mainstructure} below). 

As before smooth sections may be embedded into $\Ga_{\gs}(X,E)$ by the ``constant'' embedding
now denoted by $\Sigma$, i.e.,
$\Sigma(s)=\cl[(s)_\eps]$. Since $\CC^\infty(X)$ is a subring of $\G(X)$, 
$\Ga_{\gs}(X,E)$ can also be viewed as a $\CC^\infty(X)$-module
and the two respective module structures on the space of generalized sections are
compatible in the sense of the following commutative diagram.
\[
\begin{CD}
\CC^\infty(X)\times\Ga(X,E)     @>\sigma\times\Sigma>>       \gs(X)\times\Ga_{\gs}(X,E)  \\
@VV\cdot V     @VV\cdot V                     \\
\Ga(X,E)     @>\Sigma>>   \Ga_{\gs}(X,E) \\
\end{CD}
\]

The most important structural properties of $\G(X,E)$
are subsumed in the following results.

\btheorem  $\Ga_{\gs}(\_\,,E)$ is a fine sheaf of $\gs(\_\,)$-modules. 
\etheorem
\pr This is a straightforward generalization of the $\R^n$-case.
\ep

\btheorem\label{mainstructure} The following chain of $\CC^\infty(X)$-module
isomorphisms holds:
$$
\Ga_{\gs}(X,E) \cong \G(X)\otimes_{{\cal C}^\infty(X)}\Gamma(X,E) \cong
L_{{\cal C}^\infty(X)}\Big(\Gamma(X,E^*),\G(X)\Big)
$$
\etheorem
\pr
$\Gamma(X,E)$ is projective and finitely generated (apply \cite{GHV}, 2.23, Cor.
to each connected component),
$\Gamma(X,E^*)\cong \Gamma(X,E)^*$ (\cite{GHV}, 2.24, Rem.), and, consequently,
$\Gamma(X,E)^{**}\cong \Gamma(X,E)$ (Here $\Ga(X,E)^*$ denotes the dual
$\CC^\infty(X)$-module of $\Ga(X,E)$). Hence
$\G(X)\otimes_{{\cal C}^\infty(X)}\Gamma(X,E)$\\
$\cong L_{{\cal C}^\infty(X)}\Big(\Gamma(X,E^*),\G(X)\Big)$ follows from \cite{bourbaki-algebra},
Ch.\ II, \S 4, 2. 

Since both $\Ga_{\gs}(\_\,,E)$ and
$L_{{\cal C}^\infty(\_\,)}(\Gamma(\_\,,E^*),\G(\_\,))$ are sheaves of 
$\CC^\infty(\_)$-modules (cf.\ e.g., \cite{sheaves}, (2.2.4))
and the isomorphy of the second and third module in the above chain of course 
also holds locally,
in order to finish the proof it suffices to show that
$\G(U,E) \cong \G(U)\otimes_{{\cal C}^\infty(U)}\Gamma(U,E)$ for any
trivializing open set $U\subseteq X$. But for such a $U$ we have
$\Ga_{\gs}(U,E)\cong \G(U)^{n'}$ and $\Gamma(U,E)\cong \CC^\infty(U)^{n'}$, so the claim
follows.
\ep

\brem \label{gxerem} 
$\!\!$Endowing $\G(X)\otimes_{\CC^\infty(X)} \Ga(X,E)$ with the canonical 
$\G(X)$-module structure induced by $u_1 \cdot (u_2\otimes \xi) = (u_1u_2)\otimes \xi$,  
($u_1, u_2 \in \G(X)$, $\xi\in \Ga(X,E)$) it follows
immediately that the $\CC^\infty(X)$-module isomorphism $\Ga_{\gs}(X,E)\cong \G(X)
\otimes_{\CC^\infty(X)} \Ga(X,E)$ is in fact also a $\G(X)$-module isomorphism.
\erem

\bc \label{multlin}
Let $E_1, \dots, E_k$, $F$ be vector bundles with base manifold $X$. Then the
following isomorphism of $\CC^\infty(X)$-modules holds:
$$
\Ga_{\gs}\Big(X,L(E_1,\dots,E_k;F)\!\Big) \cong L_{\CC^\infty(X)}
\Big(\Ga(X,E_1),\dots,\Ga(X,E_k);\Ga_{\gs}(X,F)\!\Big)
$$
\ecor
\pr
By Theorem \ref{mainstructure} the right hand side can be written as
\begin{eqnarray*}
\lefteqn{
 L_{\CC^\infty(X)}\Big(\Ga(X,E_1),\dots,\Ga(X,E_k);\G(X)\otimes_{\CC^\infty(X)}\Gamma(X,F)\Big)
} \\ [.5em]
&&\cong
 L_{\CC^\infty(X)}\Big(\Ga(X,E_1)\otimes_{\CC^\infty(X)}\dots\\
&& \hphantom{mmmmmmmmm}
 \dots \otimes_{\CC^\infty(X)}\Ga(X,E_k);\G(X)\otimes_{\CC^\infty(X)}\Gamma(X,F)\Big)\\[.5em]
&&\cong
 \G(X)\otimes_{{\cal C}^\infty(X)}
 L_{\CC^\infty(X)}\Big(\Ga(X,E_1)\otimes_{\CC^\infty(X)}\dots\\
&& \hphantom{mmmmmmmmmmmmmmmm}
 \dots\otimes_{\CC^\infty(X)}\Ga(X,E_k);\Gamma(X,F)\Big)\\[.5em]
&&\cong \G(X)\otimes_{{\cal C}^\infty(X)}
 L_{\CC^\infty(X)}\Big(\Ga(X,E_1),\dots,\Ga(X,E_k);\Gamma(X,F)\Big)
\end{eqnarray*}
Here the second isomorphism holds by \cite{bourbaki-algebra}, Ch.\ II \S 4, 2.,
Prop.\ 2
since $\Ga(X,E_1)\otimes_{\CC^\infty(X)}\dots\otimes_{\CC^\infty(X)}
\Ga(X,E_k)$ is finitely generated and projective. Now 
\begin{eqnarray*}
L_{\CC^\infty(X)}\Big(\Ga(X,E_1)\dots\Ga(X,E_k);\Gamma(X,F)\Big)
\cong \Gamma\Big(X,L(E_1,\dots,E_k;F)\Big)
\end{eqnarray*}
by \cite{GHV}, 2.24, Cor.\ 2, so the claim follows from Theorem \ref{mainstructure}.
\ep

\btheorem \label{gxproj}
The $\G(X)$-module $\Ga_{\gs}(X,E)$ is finitely generated and projective.
\etheorem
\pr
Choose a vector bundle $F$ such that $E\oplus F = X\times \R^{n'}$ for some ${n'}\in \N$
(apply \cite{GHV}, 2.23 to each connected component). Then we have the following 
$\G(X)$-isomorphisms:
$$
\Ga_{\gs}(X,E)\oplus_{\G(X)} \Ga_{\gs}(X,F) \cong \Ga_{\gs}(X,X\times \R^{n'}) \cong \G(X)^{n'}
$$

It follows that the $\G(X)$-module $\Ga_{\gs}(X,E)$ is a direct summand in a finitely
generated free $\G(X)$-module, hence is projective and finitely generated
(\cite{bourbaki-algebra}, Ch.\ II, \S 2, 2., Cor.\ 1).
\ep

We will study further properties of $\Ga_{\gs}(X,E)$ as a $\G(X)$-module after 
Lemma \ref{locgstensor}.


Analogously to the earlier cases we set up coupled calculus
in order to obtain a convenient language for describing compatibility
with the distributional setting. In the following 
definition, $(.|.)$ denotes the canonical vector bundle homomorphism
\begin{eqnarray*}
 (.|.) &:=& \mbox{tr}_E\otimes \mbox{id}\\
 (E\otimes E^*)\otimes \mbox{Vol}(X) &\to& (X\times \C)\otimes \mbox{Vol}(X)
= \mbox{Vol}(X)
\end{eqnarray*}
where $\mbox{tr}_E$ is the vector bundle isomorphism induced by the pointwise
action of $v^*\in E_p^*$ on $v\in E_p$.
\bd\label{gassdef}
\begin{itemize}
\item[(i)] A generalized section $S\in\Ga_{\gs}(X,E)$ is called {\em associated to $0$},
$S\approx 0$, if for all $\mu\in\Ga_c(X,E^*\otimes\Vol(X))$ and one (hence any) 
representative $(s_\eps)_\eps$ of $S$
\[
        \lim\limits_{\eps\to 0}\int\limits_X\,(s_\eps|\mu)\,=\,0\,.
\]
\item[ (ii)]Let $S\in \Ga_{\gs}(X,E)$ and $w\in\D'(X,E)$. We say that $S$ admits 
$w$ as {\em associated 
distribution (with values in $E$)} and call $w$ the {\em distributional shadow} (or 
{\em macroscopic aspect}) of $S$ if for all 
$\mu\in\Ga_c(X,E^*\otimes\Vol(X))$ and one (hence any) representative 
\[ 
        \lim_{\eps\to 0}\int\limits_X\,(s_\eps|\mu)\,=\,w(\mu)\,,
\]
where $w(\mu)$ denotes the distributional action of $w$ on $\mu$.
In that case we use the notation $S\approx w$.
\end{itemize}
\edefi
$S\approx T:\Leftrightarrow S-T\approx 0$ defines an equivalence relation 
giving rise to a linear quotient of $\Ga_{\gs}(X,E)$. 
If $S\approx T$ we call $S$ and $T$ {\em associated 
to each other}. In complete analogy to the scalar case, 
by localization we immediately have
\bp \begin{itemize}
\item[(i)]$S\approx 0$ in $\Ga_{\gs}(X,E)$ $\Leftrightarrow$ $S^i_\al\approx 0$ in 
        $\gs(\psi_\al(V_\al))$ $\forall \al, i=1,\dots,{n'}$
\item[(ii)] $S\approx w\in \D'(X,E)$  $\Leftrightarrow$ $S^i_\al\approx w^i_\al$
        in $\gs(\psi_\al(V_\al))$ $\forall\al,i=1,\dots,{n'}$
\end{itemize}
\ep
\eprop

\bd Let $S\in\Ga_{\gs}(X,E)$.
\begin{itemize}
\item[(i)] $S$ is called $\CC^k$-associated to $0$ $(0\leq k\leq\infty)$, $S\approx_k 0$,
if for  one (hence any) representative $(s_\eps)_\eps$ and 
$\forall \al,\,i=1,\dots,{n'}$ $s^i_\al\,_\eps\to 0$ uniformly on compact sets 
in all derivatives of order less or (if $k<\infty$) equal to $k$.

\item[(ii)] We say that $S$ allows $t\in\Ga^k(X,E)$ as a $\CC^k$-associated section, 
$S\approx_k t$, if  for  one (hence any) representative $(s_\eps)_\eps$  and 
$\forall \al,\,i=1,\dots,{n'}$ $s^i_\al\,_\eps\to t^i_\al$ uniformly on compact sets 
in all derivatives of order less or (if $k<\infty$) equal to $k$.
\end{itemize}
\edefi

As is the case with $\gs(X)$ the different $\CC^\infty$-module structures of 
$\D'(X,E)$ and  $\Ga_{\gs}(X,E)$, respectively, may be reconciled at the level of association:
\bp Let $U\in\gs(X)$ and $S\in\Ga_{\gs}(X,E)$.
\begin{itemize}
\item[(i)] If $U\approx w\in\D'(X)$, $s\in\Ga(X,E)$ and either (a) $S=\Sigma(s)$
or (b) $S\approx_\infty s$, then $U\,S\approx ws$.
\item[(ii)] If $S\approx s\in\D'(X,E)$, $f\in\CC^\infty(X)$ and either (a) 
$U=\sigma(f)$ or (b) $U\approx_\infty f$, then $U\,S\approx fs$.
\item[(iii)] If $U\approx_k f$ and $S\approx_k s$ then $U\,S\approx_k fs$ 
($f\in\CC^k(X)$, $s\in\Ga^k(X,E))$.
\end{itemize}
\eprop
\pr
Simply apply Proposition \ref{consist} componentwise. \ep

\section{Generalized tensor analysis} \label{gta}

In the case where $E\to X$ is some tensor bundle $T^r_s(X)$ over the manifold $X$ we
shall use the notation ${\gs}^r_s(X)$ for $\Ga_{\gs}(X,T^r_s(X))$ and similarly for
$\Ga_{\es}$, $\Ga_{\esm}$ and $\Ga_{\ns}$. The space of smooth tensor fields will be denoted by
${\cal T}^r_s(X)$. One of the main goals in our analysis of this particular case
of generalized sections of vector bundles is to demonstrate the relative ease 
with which
arguments from classical analysis can be carried over to the generalized functions
setting. Our first result gives several algebraic characterizations 
of $\G^r_s(X)$. 

\btheorem\label{multilinear} 
\begin{itemize}
\item[(i)] As $\G(X)$-module, $\G^r_s(X)\cong
L_{\G(X)}\Big(\G^0_1(X)^r,\G^1_0(X)^s;\G(X)\Big)$. 
\item[(ii)] As $\CC^\infty(X)$-module, 
$\G^r_s(X)\cong L_{\CC^\infty(X)}\Big(\Om^1(X)^r,{\mathfrak X}(X)^s;\G(X)\Big)$. 
\item[(iii)] As $\CC^\infty(X)$-module and also as  $\G(X)$-module,\\
\[
 \G^r_s(X) \cong \G(X)\otimes_{\CC^\infty(X)}{\cal T}^r_s(X)\,.
\]
\end{itemize}
\etheorem

To simplify notations we will set $r=1=s$ in the proof. We first establish 
the following localization result.
\blem \label{locgstensor} Let $T\in L_{\gs(X)}({\gs}^0_1(X),{\gs}^1_0(X);\gs(X))$, 
$A\in{\gs}^0_1(X)$ and 
$\Xi\in{\gs}^1_0(X)$ with $\Xi|_U=0$ for some open $U\subseteq X$. Then $T(A,\Xi)|_U=0$.
\elem

\pr Since $U$ can be written as the union of a collection of open sets 
$(U_p)_{p\in U}$ such that each $\overline U_p\comp V_\al$ for some chart $V_\al$ 
and due to the sheaf property of $\gs(X)$ we may assume without loss of generality 
that $\overline U\comp V_\al$ 
and write $\Xi|_{V_\al}=\Xi^i\pa_i$ with $\Xi^i\in\gs(V_\al)$ vanishing on $U$. 
Let now $f$ be a bump function on $\overline U$ (i.e., $f\in\D(V_\al)$, 
$f|_{\overline U}=1$) then (using summation convention)
\beas
        T(A,\Xi)|_U&=&f^2|_U\,T(A,\Xi)|_U\,=\,f^2\,T(A,\Xi)|_U\\
        &=&T(A,f\Xi^if\pa_i)|_U
        \,=\,f\Xi^i\,T(A,f\pa_i)|_U\\
        &=& f\Xi^i|_U\,T(A,f\pa_i)|_U\,=\,0\,,
\eeas
where we did not distinguish notationally between $f$ and $\sigma(f)$.
\ep
\vskip6pt
From this result it follows 
that for any $V\subseteq X$ open, $A\in{\gs}^0_1(V)$ and $\Xi\in{\gs}^1_0(V)$
we may unambiguously define $T|_V(A,\Xi)$.
\vskip6pt
{\it Proof of the theorem. } (i)
Let $T=\cl[(t_\eps)_\eps]\in{\gs}^1_1(X)$, $A=\cl[(a_\eps)_\eps]\in{\gs}^0_1(X)$ and
$\Xi=\cl[(\xi_\eps)_\eps]\in{\gs}^1_0(X)$. Using classical contraction we define 
componentwise the following map
\[ \tilde T:\,(a_\eps,\xi_\eps)\mapsto f_\eps:=t_\eps(a_\eps,\xi_\eps)\,.\]
From the local description it is easy to see that $F=\cl[(f_\eps)_\eps]\in\gs(X)$,
$\tilde T:{\gs}^0_1(X)\times{\gs}^1_0(X)\to\gs(X)$ is well-defined and 
$\gs(X)$-bilinear,
so $\tilde T\in L_{\gs(X)}({\gs}^0_1(X)$, 
${\gs}^1_0(X);\gs(X))$. Moreover, the assignment
$T\mapsto\tilde T$ is also $\gs(X)$-linear, so it only remains to show that 
the latter is an isomorphism.

To prove injectivity assume $\tilde T=0$, that is $(t_\eps(a_\eps,\xi_\eps))_\eps
\in\ns(X)$ for all $A=\cl[(a_\eps)_\eps]\in{\gs}^0_1(X)$ and all $\Xi=\cl[(\xi_\eps)_\eps]
\in{\gs}^1_0(X)$. To show that $T=0\in{\gs}^1_1(X)$ it suffices to work locally. 
Choose $K\comp V_\al$ and $A\in \gs^0_1(X)$, $\Xi\in \gs^1_0(X)$ whose compact 
supports are contained in $V_\al$ and such that $A=\Sigma(dx^i)$, $\Xi=\Sigma(\pa_j)$
on an open neighborhood $U$ of $K$ in $V_\al$ ($1\leq i,j\leq n$).  
Then $\ns(U)\ni
(t_\eps(a_\eps,\xi_\eps)|_{U})_\eps=(t_\al\,^i_{j\,\eps}|_U)_\eps$. 
Since $i,j$ were arbitrary we are done.

To show surjectivity choose $\tilde T\in L_{\gs(X)}({\gs}^0_1(X),{\gs}^1_0(X);\gs(X))$.
By the remark following Lemma \ref{locgstensor}, for any chart $(V_\al,\psi_\al)$ with 
coordinates $x^i$ we may define
\[
   T_{\al}\,^i_j\,=\,\tilde T|_{V_\al}(dx^i,\pa_j)\circ\psi_{\al}^{-1}\,
\in\gs(\psi_{\al}(V_\al))\,,
\]
Since $\tilde T$ is globally defined the 
$(T_{\al })_{\al}$ form a coherent family. Hence by the sheaf property of 
${\gs}^1_1(X)$ there exists a unique $T\in{\gs}^1_1(X)$ 
represented by the family $(T_{\al })_{\al}$ and by construction $\tilde T$
is the image of $T$. 

(ii) follows from Corollary \ref{multlin} 
(alternatively, it can be proved analogously to (i)).
Finally, (iii) is immediate from Theorem \ref{mainstructure} and Remark \ref{gxerem}.
\ep

Theorem \ref{multilinear} (iii) was suggested as a {\em definition} for the space of
Colom\-be\-au tensor fields in \cite{hermannbook}, Ch.\ 2. 
The proof of Theorem \ref{multilinear} (i) is easily adapted to yield the
following result on spaces of generalized sections:
\bp \label{gsecmultlin} 
Let $E_1, \dots E_k$, $F$ be vector bundles with base manifold $X$. Then the
following isomorphism of $\G(X)$-modules holds:
$$
\Ga_{\gs}\Big(\!X,L(E_1,\dots,E_k;F)\!\Big)\!\cong\!L_{\G(X)}
\Big(\Ga_{\gs}(X,E_1),\dots,\Ga_{\gs}(X,E_k);\Ga_{\gs}(X,F)\!\Big)
$$
\ep
\eprop

(An alternative proof of Proposition \ref{gsecmultlin} can be given along the lines of
\cite{GHV}, 2.24.) Hence
\begin{equation}\label{gxdual}
L_{\G(X)}(\Ga_{\gs}(X,E),\G(X)) \cong \Ga_{\gs}(X,E^*)\,.
\end{equation}
It follows that the $\G(X)$-module $\Ga_{\gs}(X,E)$ is reflexive.
Also, we note that the proof of \cite{GHV}, Ch.\ II, Prop.\ XIV can directly
be adapted to establish:
\bp \label{gxtensorprop}
Let $E$, $F$ be vector bundles with base manifold $X$. Then the
following isomorphism of $\G(X)$-modules holds:
\begin{equation}\label{gxtensor}
\Ga_{\gs}(X,E)\otimes_{\G(X)}\Ga_{\gs}(X,F) \cong \Ga_{\gs}(X,E\otimes F)  
\end{equation}
\ep
\eprop
In particular, from (\ref{gxdual}), Proposition \ref{gxtensorprop} 
and Theorem \ref{gxproj}
we conclude:
\begin{eqnarray*}
\lefteqn{L_{\G(X)}\Big(\Ga_{\gs}(X,E_1),\dots,\Ga_{\gs}(X,E_k);\Ga_{\gs}(X,F)\Big)
}\\[.5em]
&&\cong L_{\G(X)}\Big(\Ga_{\gs}(X,E_1)\otimes_{\G(X)}\dots\\
&&\hphantom{mmmmmmmmmn}\dots\otimes_{\G(X)}\Ga_{\gs}(X,E_k);\G(X)\Big)
\otimes_{\G(X)} \Ga_{\gs}(X,F) \\[.5em]
&&\cong L_{\G(X)}\Big(L_{\G(X)}(\Ga_{\gs}(X,E_1^*)\otimes_{\G(X)}\dots\\
&&\hphantom{mmmmmm}\dots\otimes_{\G(X)}\Ga_{\gs}(X,E_k^*);
\G(X));\G(X)\Big)\otimes_{\G(X)} \Ga_{\gs}(X,F)\\[.5em]
&&\cong L_{\G(X)}\Big(\Ga_{\gs}(X,E_1),\G(X)\Big)\otimes_{\G(X)}\dots\\
&&\hphantom{mmmmnnn}\dots\otimes_{\G(X)}
L_{\G(X)}\Big(\Ga_{\gs}(X,E_k),\G(X)\Big)\otimes_{\G(X)} \Ga_{\gs}(X,F)
\end{eqnarray*}
(using \cite{bourbaki-algebra}, Ch.\ II, \S 4, 2., Prop.\ 2, 4., Prop.\ 4, Cor.\ 1,
and 2., Rem.\ (2)).

Returning now to the special case of tensor bundles, 
given a generalized tensor field $T\in{\gs}^r_s(X)$ 
we shall call the $n^{r+s}$ generalized functions on $V_\al$ defined by
\[ 
        T^\alpha\,^{i_1\dots i_r}_{j_1\dots j_s}
        \,:=\,T|_{V_\al}(dx^{i_1},\dots ,dx^{i_r},\pa_{j_1},\dots,\pa_{j_s})
\]
its {\em components} with respect to the chart $(V_\al,\psi_\al)$. We shall 
use abstract index notation (cf.\ \cite{penrose_rindler}, Chap.\ 2) 
whenever convenient and write $T^{a_1\dots a_r}_{b_1\dots b_s}\in{\gs}^r_s(X)$. 
To clearly distinguish between the notions of abstract and concrete indices we 
reserve the letters $a,b,c,d,e,f$ for the previous one and $i,j,k,l,\dots$ for the 
latter one. Hence we shall denote the
components of $\Xi^a\in{\gs}^1_0(X)$ and $A_a\in{\gs}^0_1(X)$ 
w.r.t.\ the chart $(V_\al,\psi_\al)$ by $\Xi^{\al\,i}$ and ${A^\al}_i$ 
respectively. Similarly the components of a representative 
$(t^{a_1\dots a_r}_{b_1\dots b_s}\,_\eps)_\eps\in\esm\,^r_s(X)$ of 
$T^{a_1\dots a_r}_{b_1\dots b_s}\in{\gs}^r_s(X)$ will be denoted by
$(t^\alpha\,^{i_1\dots i_r}_{j_1\dots j_s}\,_\eps)_\eps$.

The spaces of moderate respectively negligible nets of tensor fields 
may be characterized invariantly by the Lie derivative (similar to the scalar case,
cf.\ Lemma \ref{emlem} (ii)).
\bp
\beas
        \esm\,^r_s(X)\!\!&=& \!\! \{ (t_\eps)_{\eps\in I}\in\! (\es)^r_s(X):
                \forall K\subset\subset X, \,\forall k\in\!\N_0\,\exists N\in\N\,
                \,\forall \xi_1,\dots,\xi_k  \\ \nn
        && \hspace{.7cm}\in{\cal T}^1_0(X): 
           \sup_{p\in K}||L_{\xi_1}\dots L_{\xi_k}\,t_\eps(p)||=O(\eps^{-N})
                \mbox{ as }\eps\to 0\}\\[.5em]
        {\ns}^r_s(X)\!\!&=\!\!& \{ (t_\eps)_{\eps\in I}\in\! (\es)^r_s(X) : 
                \forall K\subset\subset X,\,\forall k,m \in\N_0\,
                \forall \xi_1,\dots,\xi_k\nn\\ 
        &&\hspace{.9cm} \in{\cal T}^1_0(X):
\sup_{p\in K}||L_{\xi_1}\dots L_{\xi_k}\,t_\eps(p)||=O(\eps^{m})
                \mbox{ as }\eps\to 0\}\nn
\eeas
where $||. ||$ denotes the norm induced on $\T^r_s(X)$ by any Riemannian metric
on $X$.
\eprop

\bd
Let $S\in{\gs}^r_s(X)$ and $T\in{\gs}^{r'}_{s'}(X)$. We define the {\em tensor product}
$S\otimes T\in{\gs}^{r+r'}_{s+s'}(X)$ of $S$ and $T$ by
\[ S\otimes T\,:=\,\cl[(s_\eps\otimes t_\eps)_\eps]\,.
\]
\edefi

Using the local description it is easily checked that the tensor product is well
defined. Moreover it is $\gs(X)$-bilinear, associative and by a straightforward
generalization of Proposition \ref{consist} displays the following 
consistency properties
with respect to the classical resp. distributional tensor product.
\bp Let $S\in{\gs}^r_s(X)$ and $T\in{\gs}^{r'}_{s'}(X)$.
\begin{itemize}
\item[(i)]
If $T\approx w\in\D'^{r'}_{s'}(X)$, $s\in\T^r_s(X)$ and either (a) $S=\Sigma(s)$ 
or (b) $S\approx_\infty s$  
then $S\otimes T\approx s\otimes w$ in ${\gs}^{r+r'}_{s+s'}(X)$.
\item[(ii)] If $S\approx_k s$ and $T\approx_k t$ then $S\otimes T\approx_k s\otimes t$
in ${\gs}^{r+r'}_{s+s'}(X)$ $(s\in\Ga^k(X,T^r_s(X))$, $t\in\Ga^k(X,T^{r'}_{s'}(X)))$.
\end{itemize}
\ep
\eprop
We may now easily generalize the following notions of classical tensor
calculus.
\bd\begin{itemize} 
\item[(i)]  
Let $T^{a_1\dots a_r}_{b_1\dots b_s}\in{\gs}^r_s(X)$. 
We define the contraction of $T^{a_1\dots a_r}_{b_1\dots b_s}$ 
by
\[ 
        T^{a_1\dots i\dots a_r}_{b_1\dots i\dots b_s}
        :=\cl[(t^{a_1\dots i\dots a_r}_{b_1\dots i\dots b_s}\,_\eps)_\eps]
        \in{\gs}^{r-1}_{\,s-1}(X)\,.
\]
\item[(ii)]  
For any smooth vector field $\xi$ on $X$ the Lie derivative
of $T\in{\gs}^r_s(X)$ with respect to $\xi$ is given by
\[
         L_\xi T\,:=\,\cl[(L_\xi t_\eps)_\eps]\,.
\]
\item[(iii)] 
Finally, we define the universal generalized tensor algebra over $X$ by
\[ 
        \hat G(X)\,:=\,\bigoplus\limits_{r,s}{\gs}^r_s(X)\,.
\]
\end{itemize}
\edefi 

The Lie derivative displays the following consistency property with 
respect to its distributional counterpart
\bp Let $X$ be orientable and $T\approx t$ in $\gs^r_s(X)$. Then $ L_\xi T
\approx L_\xi t$.
\ep
\eprop

Next we introduce the generalized Lie derivative, i.e., the Lie derivative 
with respect to a generalized vector field. We note that an analogous definition
(i.e., Lie derivative of a distributional tensor field with respect to a
distributional vector field)
is impossible in the purely distributional setting (cf.\ \cite{marsden}, \S 5).
\bd\label{glie}
Let $\Xi\in{\gs}^1_0(X)$ and $T\in{\gs}^r_s(X)$. We define the 
{\em generalized Lie derivative} of $T$ with respect to $\Xi$ by 
\[
        L_\Xi(T)\,:=\,\cl[(L_{\xi_\eps}(t_\eps))_\eps]\,.
\]
In case $U\in\gs(X)$ we also use the notation $\Xi(U)$ for $L_\Xi U$.
\edefi

The well-definedness of $L_\Xi(T)$ is an easy consequence of the local description.
Literally all classical (algebraic) properties of the Lie derivative carry over since 
they hold componentwise. 
In particular, for generalized vector fields $\Xi,H$ we have $L_\Xi H=[\Xi,H]:=
\cl[\,([\xi_\eps,\eta_\eps])_\eps]$ and for all generalized functions 
$U$ we have: $[U\Xi,H]=U[\Xi,H]-H(U)\Xi$. Moreover, we immediately get the following 
consistency properties.
\bp Let $\Xi\in{\gs}^1_0(X)$ and $T\in{\gs}^r_s(X)$
\begin{itemize}
\item[(i)] If $\Xi=\Sigma(\xi)$ for some $\xi\in\T^1_0(X)$ then $L_\Xi(T)=L_\xi(T)$.
\item[(ii)] If $\Xi\approx_\infty\xi\in\T^1_0(X)$ and $T\approx t\in\D'^r_s(X)$ or
conversely, if  $\Xi\approx\xi\in\D'^1_0(X)$ and $T\approx_\infty t\in\T^r_s(X)$ then 
$L_\Xi(T)\approx L_\xi t$.
\item[(iii)] If $\Xi\approx_k \xi$ and $T\approx_{k+1} t$ then  $L_\Xi(T)\approx_k 
L_\xi t$ $(\xi\in\Ga^k(X,TX)$, $t\in\Ga^{k+1}(X,T^r_s(X)))$.
\end{itemize}
\ep
\eprop

For a generalized vector field $\Xi$ the map 
$L_\Xi\equiv\Xi:\gs(X)\to\gs(X)$ is clearly 
$\R$-linear (in fact even ${\cal R}$-linear) and obeys the Leibniz rule, hence is a 
derivation on $\gs(X)$. Moreover any derivation on the algebra of generalized function 
arises this way.
\btheorem ${\gs}^1_0(X)$ is $($$\R$-linearly$)$ isomorphic to $\mbox{Der}\,(\gs(X))$.
\etheorem
\pr
It suffices to show that for any derivation $\theta$ on $\gs(X)$
we may construct a unique generalized vector field $\Xi$ such that $\theta(U)=\Xi(U)$
for all $U\in\gs(X)$.
We start by showing that $\theta$ is a local operator, i.e., that $U=0$ on 
$V(\subseteq X)$ open implies $\theta(U)|_V=0$. 
To this end choose any open $W$ with
$\overline{W}\comp V$  and a function $f\in\D(V)$ equal to $1$ on 
$W$. Then $U=(1-f)U$ and
\[
        \theta(U)|_W\,=\,\theta(1-f)U|_W+(1-f)\theta(U)|_W\,=\,0\,\in\gs(W)
\]
Since $\gs$ is a sheaf, $\theta(U)|_V=0$.
Now let $(V_\al,\psi_\al)$ be a chart in $X$, $x=\psi_\al(p)$ and $U\in\gs(X)$. Then
for $y$ in a neighborhood of $x$
\beas 
        &&(U\circ\psi_\al^{-1})(y)
        =(U\circ\psi_\al^{-1})(x)
        +\int\limits_0^1\,\frac{d}{dt}(U\circ\psi_\al^{-1})(x+t(y-x))\,dt\nn\\
        &&=(U\circ\psi_\al^{-1})(x)
        +\sum\limits_{i=1}^n(y^i-x^i)
        \int\limits_0^1\,D_i(U\circ\psi_\al^{-1})(x+t(y-x))\,dt\,.
\eeas
Hence in a neighborhood of $p$ ($q=\psi_\al^{-1}(y)$),
$
        U(q)\,=\,U(p)+\sum\limits_{i=1}^n(\psi_\al^i(q)-\psi_\al^i(p))\,g_i(q)\,,
$
where $g_i$ is given by the integral above whence, in particular, 
$g_i(p)=\frac{\pa}{\pa x^i}(U\circ\psi_\al^{-1})|_{x}$. Consequently
\[
        (\theta(U))(p)\,=\,\sum\limits_{i=1}^n
        \pa_i U(p)\,\theta(\psi_\al^i)(p)
\]
and we define $\Xi$ locally to be given by $\Xi^i_\al=\theta(\psi_\al^i)$
(this is well-defined by the first part of the proof).
It is easily checked that this indeed defines a coherent family in the sense
of~(\ref{gloc}).
\ep

\section{Exterior Algebra, Hamiltonian Mechanics} \label{extalg}

In this section we are going to study generalized sections 
of the bundle $\Lambda^kT^*X$, i.e., generalized $k$-forms, thereby setting
the stage for nonsmooth Hamiltonian mechanics.

To simplify notations we set $\bigwedge^k_{\gs}(X):=\Ga_{\gs}(X,\Lambda^kT^*X)$
and similar for the spaces of moderate resp.\ negligible nets of $k$-forms.
If $X$ is oriented (with its orientation induced by $\theta$) it follows from
the local description of generalized sections that $\Sigma(\om)\approx j(\om)$ for 
all $\om\in\Om^k(X)$,
where $j$ is the embedding of regular objects into the space of distributional 
$k$-forms from \cite{marsden} (see (\ref{marsdenemb})). The basic operations of
exterior algebra are carried over to our setting by componentwise 
definitions.
\bd Let $A=\cl[(\al_\eps)_\eps]\in\bigwedge^k_{\gs}(X)$, $B=\cl[(\beta_\eps)_\eps]\in
\bigwedge^l_\gs(X)$ and $\Xi=\cl[(\xi_\eps)_\eps]\in\gs^1_0(X)$. We define the
exterior derivative, the wedge product and the insertion operator, respectively, by:
\begin{itemize}
\item [(i)] $dA:=\cl[(d\al_\eps)_\eps]\in\bigwedge_\gs^{k+1}(X)$
\item [(ii)] $A\wedge B:=\cl[(\al_\eps\wedge\beta_\eps)_\eps]\in\bigwedge_\gs^{k+l}(X)$
\item [(iii)] $i_\Xi A:=\cl[(i_{\xi_\eps}\al_\eps)_\eps]\in\bigwedge_\gs^{k-1}(X)$
\end{itemize}
\edefi

Of course all the classical relations remain valid in our framework where
(in contrast to the distributional setting) in every multilinear operation all factors 
may be generalized; in particular for $A\in\bigwedge^k_{\gs}(X)$ and $\Xi,\Xi_1,\dots
\Xi_k\in\gs^1_0(X)$
we have $(\iota_\Xi A)(\Xi_2,\dots,\Xi_k)=A(\Xi,\Xi_2,\dots,\Xi_k)$ and
$L_\Xi=d\circ i_\Xi+i_\Xi\circ d$. 

A generalized $k$-form $A$ is called closed if $dA=0$ and exact if there exists
$B\in\bigwedge_\gs^{k-1}(X)$ with $dB=A$. Clearly every exact generalized $k$-form
is closed. The converse---as in the smooth case---holds locally:
\btheorem (Poincar\'e Lemma)\\
Let $A\in\bigwedge^k_{\gs}(X)$ closed. Then for each $p\in X$ there exists a neighborhood
$U$ of $p$ and $B\in\bigwedge_\gs^{k-1}(X)$ such that
\[ A|_U=dB|_U\,. \] 
\etheorem 
\pr Since it suffices to work in a local chart we may suppose that 
$U\subseteq\R^n$ is a ball around zero. 
Let $(\al_\eps)_\eps$ denote a representative of $A$. Then $d\al_\eps=n_\eps\in
\bigwedge_\ns^{k+1}(U)$. Analogous to the classical proof (cf.\ e.g., \cite{AM}, 2.4.17)
we define an operator $H: \Om^k(U) \to \Om^{k-1}(U)$ by
\[ 
 H\om(x)(v_1,\dots,v_{k-1})=\int\limits_0^1 t^{k-1}\om(tx)
 (x,v_1,\dots,v_{k-1})\,dt,
\]
where $v_1,\dots,v_{k-1}\in\R^n$. 
Then $d\circ H+H\circ d=\id$, so $\al_\eps=H d \al_\eps+
dH\al_\eps$ for each $\eps>0$. It is immediate from the explicit form of $H$ that 
$H(\bigwedge_{\esm}^k)(U)\subseteq \bigwedge_{\esm}^{k-1}(U)$ and
$H(\bigwedge_\ns^{k}(U))\subseteq \bigwedge_\ns^{k-1}(U)$.
Thus $H\al_\eps \in \bigwedge_{\esm}^{k-1}(U)$,
$H(d\al_\eps)\in \bigwedge_\ns^{k}(U)$ and, consequently,
$A = d(H A)$ in $\bigwedge_\gs^{k}(U)$.
\ep

In what follows we suppose $X$ to be oriented. 
Analogous to Definition \ref{densityint}, 
for $K\subset\subset X$, $A\in\bigwedge_\gs^n(X)$
we define the integral of $A:=\cl[(\al_\eps)_\eps]$ over $K$ by
\[\int\limits_KA:=\cl[(\,\int\limits_K \al_\eps)_\eps].\]
For $A$ compactly supported we set $\int_XA=
\int_LA$ where $L$ is any compact neighborhood of $\supp(A)$. This notion
of integration is compatible with the one introduced by Marsden for compactly supported
distributional $n$-forms (cf.\ \cite{marsden}, 2.6). More precisely, let $\al\in
\Om^n_c(X)'$ and $A\approx \al$. Then $\int A\approx\int\al$.

Also, Stokes' theorem is easily generalized to the new setting by component\-wise
application of the classical theorem.
\btheorem Let $X$ be a manifold with boundary and $A\in\bigwedge_\gs^{n-1}(X)$ with
compact support. Then
\[ \int\limits_X dA=\int\limits_{\partial X} A \]
\ep
\etheorem

Let us now turn to the task of generalizing symplectic geometry.
Let $(X,\om)$ be a symplectic manifold, i.e., suppose that $X$ is furnished
with smooth
nodegenerate and closed $2$-form $\om$. Generalizing
$\om$ to be distributional or even an element of $\bigwedge_\gs^2$
does not seem feasible since in that setting a distributional analogue of Darboux'
theorem is not attainable (cf.\ \cite{marsden}, \S 7). 
However, by Theorem \ref{multilinear}, 
$\om\in\Om^2(X)\subseteq\bigwedge_\gs^2(X)$ induces a $\gs(X)$-bilinear 
alternating map $\gs^1_0(X)\times\gs^1_0(X)\to\gs(X)$. This in turn
allows us to define the following extension of the classical isomorphism
between vector fields and one-forms induced by $\om$:
\beas
  \om_\flat:\gs^1_0(X)&\to&\gs^0_1(X)\\
  \om_\flat(\Xi)(H)&:=&\om(\Xi,H).
\eeas

This map is even a $\gs(X)$-linear isomorphism. We denote its inverse
by $\om_\sharp$ and set $\Xi^\flat=\om_\flat(\Xi)$ and $A^\sharp=\om
_\sharp(A)$. Then we have $\Xi^\flat=i_\Xi\om\in\gs^0_1(X)$, 
$\Xi^\flat(Z)=-\Xi(Z^\flat)\in\gs(X)$ and $A^\sharp(B)=-A(B^\sharp)\in\gs(X)$
for $A,B\in\gs^0_1(X)$ and $\Xi,Z\in\gs^1_0(X)$. Moreover, if $\Xi\approx\xi\in
\D^1_0(X)$ resp.\ $A\approx\al\in\D^0_1(X)$ then $\Xi^\flat\approx\xi^\flat$ resp.\
$A^\sharp\approx\al^\sharp$.

For any $H\in\gs(X)$ we call the generalized vector field
defined by
\[ \Xi_H:=(dH)^\sharp\]
the generalized Hamiltonian vector field with energy function $H$. If 
$H\approx h\in\D'(X)$ then we have $\Xi_H\approx X_h$, where $X_h$ is defined
according to \cite{marsden}, Prop. 7.3. 

Let $F=\cl[(f_\eps)_\eps]$, $G=\cl[(g_\eps)_\eps]\in\gs(X)$. We define the Poisson 
bracket of $F$ and $G$ by \[\{F,G\}:=\cl[(\{f_\eps,g_\eps\})_\eps].\]
Literally all classical properties carry over. In particular, $\{\, ,\, \}$ is
antisymmetric, the Jacobi identity holds and we have 
$\{F,G\}=L_{\Xi_G}F=-L_{\Xi_F}G=-i_{\Xi_F}i_{\Xi_G}\om$ and
$\Xi_{\{F,G\}}=-[\Xi_F,\Xi_G]$. We note that in contrast to the distributional
setting (\cite{marsden}, Prop.\ 7.4), where ill-defined products of distributions
have to be avoided carefully, in our present framework
{\em both} factors $F$ and $G$ may be generalized functions. There is of course
a result analogous to Proposition \ref{consist} concerning consistency with 
respect to the 
smooth resp.\ distributional setting in the sense of association. 

\bex
We close this section by discussing a simple example from nonsmooth mechanics to
indicate the usefulness of the present setting. Let $X=\R^2$ and consider the
generalized Hamiltonian function $H(p,q)=\frac{p^2}{2}+D(q)$, where $D$ denotes
a generalized delta function in the sense of \cite{ode_prep}, i.e., we suppose
that $D$ possesses a 
representative $\de_\eps$ with $\supp(\de_\eps)\to\{0\}$, $\int\de_\eps\to 1$ 
and $\int|\de_\eps|\leq C$ for $\eps$ small. Clearly, every generalized delta
function is associated to $\delta$. Nets $\delta_\eps$ possessing the above
mentioned properties provide a general and flexible means of modeling delta-type
singularities (so-called {\em strict delta nets}, cf.\ \cite{MObook}, chap.\ II,
\S 7). The Hamiltonian equations for this setup take the form
$$
\dot p=-\frac{\pa H}{\pa q}=-D'(q) \qquad \dot q=\frac{\pa H}{\pa p}=p\,,
$$
leading to 
\beq\label{ode} 
\begin{array}{c}
\ddot q+D'[q]=0\\
q(0)=q_0,\, \dot q(0) = \dot q_0.
\end{array}
\eeq
This initial value problem
has been studied in detail in \cite{ode_prep,HO}. It was
shown that provided $D$ satisfies certain growth restrictions,
a solution in the Colombeau algebra exists and is unique for arbitrary initial conditions 
$q_0,\, \dot q_0 \in {\cal R}$. The limiting behavior of this unique solution
will in general depend on the chosen regularization for $\delta$. 
For example, if we choose 
$\de_\eps(x)=\frac{1}{\eps}\rho(\frac{x}{\eps})$ with $\rho\in\D(\R)$ 
we get the picture of pure reflection at the origin, i.e., 
the unique solution to (\ref{ode}) is associated to the function 
$t\to \mbox{sign}(q_0)|q_0 + \dot q_0 t|$.
(The proof consists in a
rather technical analysis of the limiting behavior of the trajectories, 
establishing that they are neither delayed nor trapped at the origin
as $\eps\to 0$.) For generalized delta functions of different type,
a more complicated limiting behavior can be observed:
For any given finite subset $S$ of $(0,\infty)$ there exists a generalized
delta function such that the solutions to (\ref{ode}) with $x_0\not=0$
and $\dot x_0 = -\mbox{sign}(x_0)\sqrt{2s}$ with $s\in S$
are trapped at the origin after time $t = -\frac{x_0}{\dot x_0}$.

Furthermore, (\ref{ode}) possesses a unique flow which itself is a Colombeau
generalized function.
Although problematic in the distributional picture (\cite{marsden}, \S 8),
energy conservation in our present setting is immediate from $\{H,H\}=0$.
\eex

The main applications of Colombeau's special algebra on manifolds 
so far have occurred in general relativity with the purpose of 
studying singular spacetimes (see \cite{vickersESI} for a survey). 
Based on the framework developed in the
present article, a satisfying theory for analyzing the geometry of these spacetimes
can be given. A thorough investigation of such generalized semi-Riemannian
geometries is deferred to a separate paper (\cite{curvature}).
\vskip12pt

\noindent{\bf Acknowledgements:}
We are indebted to M. Oberguggenberger and M. Grosser for many fruitful discussions
and their ongoing support and to A. Cap and A. Kriegl for helpful comments.
This work was in part supported by research grant P12023-MAT of the Austrian Science 
Foundation (FWF).

\end{document}